\documentclass[preprint,authoryear,12pt]{elsarticle}

\usepackage{graphicx}
\usepackage{amssymb}
\usepackage{amsmath,amsfonts,amsbsy}

\journal{Computational Statistics \& Data Analysis}

\begin{document}

\begin{frontmatter}

\title{Robust test statistics for the two-way MANOVA based on
the minimum covariance determinant estimator}

\author[label1]{Bernhard Spangl\corref{cor1}}

\ead{bernhard.spangl@boku.ac.at}

\cortext[cor1]{Corresponding author.}

\address[label1]{Institute of Applied Statistics and Computing,
  University of Natural Resources and Life Sciences, Vienna,
  Gregor-Mendel-Str.~33, 1180 Vienna, Austria}

\begin{abstract}

  Robust test statistics for the two-way MANOVA based on the minimum
  covariance determinant (MCD) estimator are proposed as alternatives
  to the classical Wilks' Lambda test statistics which are well known
  to be very sensitive to outliers as they are based on classical
  normal theory estimates of generalized variances. The classical
  Wilks' Lambda statistics are robustified by replacing the
  classical estimates by highly robust and efficient reweighted MCD
  estimates.
  Further, Monte Carlo simulations are used to evaluate the performance
  of the new test statistics under various designs by investigating
  their finite sample accuracy, power, and 
  robustness against outliers. 
  Finally, these robust test statistics are applied to a real data
  example.

\end{abstract}


\end{frontmatter}

\section{Introduction}
\label{sec:intro}

Two-way multivariate analysis of variance (MANOVA) deals with testing
the effects of the two grouping variables, usually called factors,
on the measured observations as well as interaction effects between
the factors. 
It is the direct multivariate analogue of two-way univariate
ANOVA and is able to deal with possible correlations between
the variables under consideration. 
The classical MANOVA is based on the decomposition
of the total sum of squares and cross-products (SSP) matrix,  
and test statistics
then particularly compare some scale measures between two
matrices, such as the determinant, the trace, or the largest
eigenvalue of the matrices.

We now consider the two-way fixed-effects layout.
Let us suppose we have $N = rcn$
independent $p$-dimensional observations generated by the model
\begin{eqnarray}
  \label{eq:MANOVAwithInter}
  \boldsymbol y_{ijk} = \boldsymbol \mu + 
  \boldsymbol \alpha_i + \boldsymbol \beta_j +
  \boldsymbol \gamma_{ij} + \boldsymbol e_{ijk} \ ,
\end{eqnarray}
with $i = 1, \ldots, r$, $j = 1, \ldots, c$, $k = 1, \ldots, n$,
where $\boldsymbol \mu$ is an overall effect, 
$\boldsymbol \alpha_i$ is the $i$-th level effect of the first factor with
$r$ levels,
and $\boldsymbol \beta_j$ is the $j$-th level effect of the second
factor with $c$ levels.
We will call the first and second factor row and column factor,
respectively. 
$\boldsymbol \gamma_{ij}$ is the interaction
effect between the $i$-th row factor level and the $j$-th column factor
level, and $\boldsymbol e_{ijk}$ is the error term which is assumed
under classical assumptions 
to be independent and identically distributed as
$\mathcal N_p (\boldsymbol 0, \boldsymbol \Sigma)$ for all $i$, $j$, $k$.
Additionally, we have $\sum_{i=1}^r \boldsymbol \alpha_i = \sum_{j=1}^c
\boldsymbol \beta_j = \sum_{i=1}^r \boldsymbol \gamma_{ij} = \sum_{j=1}^c
\boldsymbol \gamma_{ij} = \boldsymbol 0$. We require that the number of
observations in each factor combination group should be the same, so that
the total SSP matrix can be suitably decomposed.

We are interested in testing the null hypotheses of all
$\boldsymbol \alpha_i$ being equal to zero,
all $\boldsymbol \beta_j$ being equal to zero, and
all $\boldsymbol \gamma_{ij}$ being equal to zero.
One of the most widely used tests in
literature to decide upon these hypotheses 
is the likelihood ratio (LR) test. Here, the LR test statistic is
distributed according to a Wilks' Lambda distribution
$\mathit \Lambda (p, \nu_1, \nu_2)$ with
parameters $p$, $\nu_1$, and $\nu_2$ \citep{WILKS1932}. We will give
the details in the following and start with the test for interactions.
The hypothesis we want to test is
\begin{eqnarray}
  \label{eq:HAB}
  H_{AB} : \boldsymbol \gamma_{11} = \boldsymbol \gamma_{12} = \ldots =
  \boldsymbol \gamma_{rc} = \boldsymbol 0 
\end{eqnarray}
against the alternative that at least one $\boldsymbol \gamma_{ij} \ne
\boldsymbol 0$. 
Further, we may also test for main effects. 
E.g., the hypothesis
we want to test is that there are no row effects, i.e.,
\begin{eqnarray}
  \label{eq:HA}
  H_{A} : \boldsymbol \alpha_{1} = \boldsymbol \alpha_{2} = \ldots =
  \boldsymbol \alpha_{r} = \boldsymbol 0 
\end{eqnarray}
against the alternative that at least one $\boldsymbol \alpha_{i} \ne
\boldsymbol 0$. 
A similar hypothesis may be tested for column effects $\boldsymbol \beta_j$.
More details may be found in, e.g., \citet{MARDIA1979}.

It is widely known that test statistics based on sample covariances
such as Wilks' Lambda 
are very sensitive to outliers. Therefore, inference based on such
statistics can be severely distorted when the data is contaminated
by outliers. 
A common approach to robustify statistical inference procedures is
to replace the classical non-robust estimates by robust ones. The
underlying idea of this plug-in principle is that robust
estimators reliably estimate the parameters of the distribution of
the majority of the data and that this majority follows the
classical model.
Hence, robust test statistics require robust estimators of
scatter instead of sample covariance matrices. Many such robust
estimators of scatter have been proposed in the literature. 
See, e.g., \citet{HUBERT2008} for an overview. 

In this article we use highly robust and efficient reweighted 
minimum covariance determinant (MCD) estimates 
to replace the classical non-robust
covariance estimates of the Wilks' Lambda test statistic
and, in this way, 
extend the approach of \citet{TODOROV2010} to two-way MANOVA designs.

Further 
approaches to robustify one-way MANOVA
have already been proposed in the literature 
to overcome the non-robustness of the classical Wilks' Lambda test
statistic in the case of data contamination.
Instead of using the original measurements 
\citet{NATH1985} suggested 
to apply the classical Wilks' Lambda statistic to the ranks of the
observations.
Another approach proposed by \citet{VANAELST2011} uses multivariate
S-estimators and the related MM-estimators. The null distribution
of the corresponding likelihood ratio type test statistic is obtained
by a fast and robust bootstrap (FRB) method. 
Whereas 
\citet{WILCOX2012} suggested a one-way MANOVA testing procedure
based on trimmed means and Winzorized covariance matrices.
\citet{XU2008} proposed
Wilks' Lambda testing procedures for one-way and
two-way MANOVA designs 
based on estimators of the covariance
matrix where the mean is replaced by the median, taken component-wise
over the indices. Although the vector of marginal medians is easily
computed it is not affine equivariant. 

The rest of this article is organized as follows. 
In Section \ref{sec:robustWilksLambda} we propose a
robust version of the Wilks' Lambda test statistic based on the 
reweighted MCD estimator.
As the distribution of this robust Wilks' Lambda test statistic
differs from the classical one, it is necessary to find a good
approximation. 
Section \ref{sec:approxDist} summarizes the construction of an approximate
distribution based on Monte Carlo simulation.
In Section \ref{sec:simulations} the design of the simulation study
is described and
Section \ref{sec:results} investigates the finite sample robustness
and power of the proposed tests and compare their performance
with that of the classical and rank transformed Wilks' Lambda tests.
In Section \ref{sec:example} the proposed tests are applied to a real
data example and
Section \ref{sec:conclusions} finishes with a discussion and directions for
further research.
Additional numerical results of the simulation study are given in
\ref{sec:appendix}.

\section{The robust Wilks' Lambda statistics}
\label{sec:robustWilksLambda}

It is well known that the Wilks' Lambda statistic---as it is based on
SSP matrices---is prone to outliers.
In order to obtain a robust procedure with high breakdown point in
the two-way MANOVA model we construct a robust version of the Wilks'
Lambda statistic by replacing the classical estimates of all SSP
matrices by reweighted MCD estimates.
The MCD estimator introduced by
\citet{ROUSSEEUW1985} looks for a subset of $h$ observations with lowest
determinant of the sample covariance matrix. The size of the subset
$h$ defines the so called trimming proportion and it is chosen between
half and the full sample size depending on the desired robustness
and expected number of outliers. The
MCD location estimate is defined as the mean of that subset and the
MCD scatter estimate is a multiple of its covariance matrix. The
multiplication factor is selected so that it is consistent at the
multivariate normal model and unbiased at small samples
\citep[cf.][]{PISON2002}.
This estimator is not very efficient at normal models,
especially if $h$ is selected so that maximal breakdown point is achieved,
but in spite of its low efficiency it is the most widely used robust
estimator in practice. To overcome the low efficiency of the MCD
estimator, a reweighed version is used.
This approach is along the same lines as proposed by \citet{TODOROV2010}
for the one-way MANOVA model.

We start by computing initial estimates of the factor combination group means
$\widehat{\boldsymbol \mu}_{ij\boldsymbol{\cdot}}^0$, 
$i = 1, \ldots, r$, $j = 1, \ldots, c$, 
and the common covariance matrix, $\boldsymbol C_0$.
The means $\widehat{\boldsymbol \mu}_{ij\boldsymbol{\cdot}}^0$
are the robustly estimated
centers of the factor combination groups based on the reweighted MCD
estimator.
To obtain $\boldsymbol C_0$
we first center each observation $\boldsymbol y_{ijk}$ by its
corresponding factor combination group mean
$\widehat{\boldsymbol \mu}_{ij\boldsymbol{\cdot}}^0$.
Then, these centered observations are pooled to get a robust
estimate of the common covariance matrix, 
$\boldsymbol C_0$, by again using the reweighted MCD estimator. 

Using the obtained estimates
$\widehat{\boldsymbol \mu}_{ij\boldsymbol{\cdot}}^0$ and
$\boldsymbol C_0$ we can now calculate the initial robust distances
\citep{ROUSSEEUW1991}
\begin{eqnarray}
  \label{eq:mahalanobisDistance}
  \mathit \Delta_{ijk}^0 = \sqrt{
  (\boldsymbol y_{ijk} - \widehat{\boldsymbol \mu}_{ij\boldsymbol{\cdot}}^0)^\top
  \boldsymbol C_0^{-1}
  (\boldsymbol y_{ijk} - \widehat{\boldsymbol \mu}_{ij\boldsymbol{\cdot}}^0) } \ .
\end{eqnarray}
With these initial robust distances we are able to define a weight for each
observation $\boldsymbol y_{ijk}$, $i=1, \ldots, r$, $j=1, \dots, r$,
$k=1, \ldots, n$ by setting the weight to 1 if the corresponding robust
distance is less or equal to $\sqrt{ \chi_{p ; 0.975}^2 }$ and to 0
otherwise, i.e.,
\begin{eqnarray}
  \label{eq:weights}
  w_{ijk} = \left \{
  \begin{array}{ll}
    1 & \mbox{if } \mathit \Delta_{ijk}^0 \le \sqrt{ \chi_{p ; 0.975}^2 } \\
    0 & \mbox{otherwise} \ .
  \end{array} \right.
\end{eqnarray}
With these weights we can calculate the final reweighted estimates
which are necessary for constructing the robust Wilks' Lambda
statistics $\mathit \Lambda_R^{\boldsymbol \cdot}$:
\begin{eqnarray*}
  \widehat{\boldsymbol \mu}_{ij\boldsymbol{\cdot}} & = & 
  \frac{1}{w_{ij\boldsymbol{\cdot}}} \sum_{k=1}^n w_{ijk} \boldsymbol y_{ijk} \ , \quad
  \widehat{\boldsymbol \mu}_{i\boldsymbol{\cdot \cdot}} =
  \frac{1}{w_{i\boldsymbol{\cdot\cdot}}} \sum_{j=1}^c \sum_{k=1}^n w_{ijk} \boldsymbol y_{ijk} \ , \quad
  \widehat{\boldsymbol \mu}_{\boldsymbol{\cdot} j \boldsymbol{\cdot}} =
  \frac{1}{w_{\boldsymbol{\cdot}j\boldsymbol{\cdot}}} \sum_{i=1}^r \sum_{k=1}^n w_{ijk} \boldsymbol y_{ijk} \ , \\
  \widehat{\boldsymbol \mu}_{\boldsymbol{\cdot \cdot \cdot}} & = & 
  \frac{1}{w_{\boldsymbol{\cdot\cdot\cdot}}} \sum_{i=1}^r \sum_{j=1}^c \sum_{k=1}^n w_{ijk} \boldsymbol y_{ijk} \ ,\\
  \boldsymbol W_R & = & \sum_{i=1}^r \sum_{j=1}^c \sum_{k=1}^n
    w_{ijk} (\boldsymbol y_{ijk} -
     \widehat{\boldsymbol \mu}_{ij \boldsymbol{\cdot}})
    (\boldsymbol y_{ijk} -
     \widehat{\boldsymbol \mu}_{ij \boldsymbol{\cdot}})^\top \ ,\\
  \boldsymbol E_R & = & \sum_{i=1}^r \sum_{j=1}^c \sum_{k=1}^n
    w_{ijk} (\boldsymbol y_{ijk} -
     \widehat{\boldsymbol \mu}_{i\boldsymbol{\cdot \cdot}} -
     \widehat{\boldsymbol \mu}_{\boldsymbol{\cdot} j \boldsymbol{\cdot}} +
     \widehat{\boldsymbol \mu}_{\boldsymbol{\cdot \cdot \cdot}})
    (\boldsymbol y_{ijk} -
     \widehat{\boldsymbol \mu}_{i\boldsymbol{\cdot \cdot}} -
     \widehat{\boldsymbol \mu}_{\boldsymbol{\cdot} j \boldsymbol{\cdot}} +
     \widehat{\boldsymbol \mu}_{\boldsymbol{\cdot \cdot \cdot}})^\top \ ,
                        \mbox{and} \\
  \boldsymbol R_R & = & \sum_{i=1}^r
     w_{i\boldsymbol{\cdot\cdot}}
     (\widehat{\boldsymbol \mu}_{i\boldsymbol{\cdot \cdot}} -
     \widehat{\boldsymbol \mu}_{\boldsymbol{\cdot \cdot \cdot}})
    (\widehat{\boldsymbol \mu}_{i\boldsymbol{\cdot \cdot}} -
     \widehat{\boldsymbol \mu}_{\boldsymbol{\cdot \cdot \cdot}})^\top \ ,
\end{eqnarray*}
where
\begin{eqnarray*}
  w_{ij\boldsymbol{\cdot}} =  
  \sum_{k=1}^n w_{ijk} \ , \quad
  w_{i\boldsymbol{\cdot \cdot}} =
  \sum_{j=1}^c \sum_{k=1}^n w_{ijk} \ , \quad
  w_{\boldsymbol{\cdot} j \boldsymbol{\cdot}} =
  \sum_{i=1}^r \sum_{k=1}^n w_{ijk} \ , \\
\end{eqnarray*}
and
\begin{eqnarray*}
  w_{\boldsymbol{\cdot \cdot \cdot}} = 
  \sum_{i=1}^r \sum_{j=1}^c \sum_{k=1}^n w_{ijk} \ .
\end{eqnarray*}
The robust test statistic for testing the equality of the interaction
terms, cf.~(\ref{eq:HAB}), is
\begin{eqnarray}
  \label{eq:HABstatistic}
  \mathit \Lambda_R^{AB} = 
  \frac{| \boldsymbol W_R |}{| \boldsymbol E_R |} \ , 
\end{eqnarray}
where $| \boldsymbol S |$ is the determinant of $\boldsymbol S$.
We reject the hypothesis of no interactions for small values of
$\mathit \Lambda_R^{AB}$.
Further,
the robust test statistic to test the equality of the $\boldsymbol \alpha_{i}$
being zero 
irrespective of the $\boldsymbol \beta_{j}$ and $\boldsymbol \gamma_{ij}$, 
cf.~(\ref{eq:HA}), is
\begin{eqnarray}
  \label{eq:HAstatisticWith}
  \mathit \Lambda_R^{A} = 
  \frac{| \boldsymbol W_R |}{| \boldsymbol W_R + \boldsymbol R_R |} \ . 
\end{eqnarray}
Equality is again rejected for small values of $\mathit \Lambda_R^A$.
We note that if significant interactions exist then it does not make 
much sense to test for row and column effects. 

Alternatively, we might decide to ignore interaction effects completely.
In such a case, the two-way MANOVA design without interactions,
\begin{eqnarray}
  \label{eq:MANOVAwithoutInter}
  \boldsymbol y_{ijk} = \boldsymbol \mu + 
  \boldsymbol \alpha_i + \boldsymbol \beta_j +
  + \boldsymbol e_{ijk} \ ,
\end{eqnarray}
we consider the
error matrix $\boldsymbol E_R$ instead of $\boldsymbol W_R$ and the
robust test statistic for testing for row effects is then
\begin{eqnarray}
  \label{eq:HAstatisticWithout}
  \mathit \Lambda_R^{A} = 
  \frac{| \boldsymbol E_R |}{| \boldsymbol E_R + \boldsymbol R_R |} \ .
\end{eqnarray}

In similar manner we may define robust tests for column effects
$\boldsymbol \beta_j$. Setting
all weights equal to 1, the robust test statistics coincide with the
classic ones. 

For computing the MCD and related estimators the FAST-MCD algorithm
of \citet{ROUSSEEUW1999} will be used as implemented in the R package
\texttt{rrcov} \citep[cf.][]{TODOROV2009}.

\section{Approximate distribution of
  $\mathit \Lambda_R^{\boldsymbol \cdot}$}
\label{sec:approxDist}

The Wilks' Lambda distribution 
$\mathit \Lambda (p, \nu_1, \nu_2)$ 
was originally proposed by \citet{WILKS1932}.
For some special cases, functions of statistics
that have a Wilks' Lambda distribution
may be expressed in terms of the $F$~distribution
\citep[cf.][pp.~195-196]{ANDERSON1958}.
For other cases, provided $\nu_1$ is large, one of the most popular
approximation is Bartlett's $\chi^2$~approximation \citep{BARTLETT1938}
given by
\begin{eqnarray}
  \label{eq:BartlettApprox}
  - \left ( \nu_1 - \frac{1}{2} (p - \nu_2 + 1) \right ) \ln
  \mathit \Lambda (p, \nu_1, \nu_2) \approx \chi^2_{p \nu_2} \ .
\end{eqnarray}
Analogously to this $\chi^2$~approximation of the classical
statistic and as proposed by \citet{TODOROV2010} for the one-way MANOVA
we will use the following approximation for
$\mathit \Lambda_R^{\boldsymbol \cdot}$:
\begin{eqnarray}
  \label{eq:robustApprox}  
  L_R^{\boldsymbol \cdot} =
  - \ln \mathit \Lambda_R^{\boldsymbol \cdot} \approx \delta \chi_q^2 \ ,
\end{eqnarray}
and express the multiplication factor $\delta$ and the degrees of freedom
of the $\chi^2$~distribution, $q$, through the expectation and variance
of $L_R^{\boldsymbol \cdot}$, i.e.,
\begin{eqnarray}
\label{eq:robustApproxExpectVar}  
  \mathbb{E} (L_R^{\boldsymbol \cdot}) = \delta q \quad \mbox{and} \quad
  \mbox{Var} (L_R^{\boldsymbol \cdot}) = 2 \delta^2 q \ .
\end{eqnarray}
Hence, we get
\begin{eqnarray}
\label{eq:robustApproxDQ}  
  \delta = \mathbb{E} (L_R^{\boldsymbol \cdot}) \frac{1}{q} \quad \mbox{and} \quad
  q = 2 \frac{\mathbb{E} (L_R^{\boldsymbol \cdot})^2}{\mbox{Var}
  (L_R^{\boldsymbol \cdot})} \ .
\end{eqnarray}
$\mathbb{E} (L_R^{\boldsymbol \cdot})$ and $\mbox{Var} (L_R^{\boldsymbol \cdot})$
are determined by simulation as follows.

For a given dimension $p$, number of levels $r$, $c$, and sample size $n$
in each factor combination group, samples $\mathcal Y^{(\ell)} = \{
\boldsymbol y_{111}, \dots, \boldsymbol y_{ijk}, \ldots, \boldsymbol y_{rcn} \}$
of size $N = rcn$ from the $p$-variate standard normal distribution will
be generated, i.e., $\boldsymbol y_{ijk} \sim
\mathcal N_p (\boldsymbol 0, \boldsymbol I_p)$, $i = 1, \ldots, r$,
$j = 1, \ldots, c$, $k = 1, \ldots, n$. For each sample the robust
Wilks' Lambda statistic $\mathit \Lambda_R^{\boldsymbol \cdot}$ based on
the weighted MCD will be calculated. After performing $m' = 3000$ trials,
the sample mean and variance of $L_R^{\boldsymbol \cdot}$ will be obtained
as
\begin{eqnarray*}
  \mbox{ave}(L_R^{\boldsymbol \cdot}) & = & \frac{1}{m'}
  \sum_{\ell=1}^{m'} L_R^{\boldsymbol \cdot (\ell)} \quad \mbox{and} \\
  \mbox{var}(L_R^{\boldsymbol \cdot}) & = & \frac{1}{m' - 1}
  \sum_{\ell=1}^{m'} \left ( L_R^{\boldsymbol \cdot (\ell)} -
  \mbox{ave}(L_R^{\boldsymbol \cdot}) \right )^2 \ .
\end{eqnarray*}
Substituting these values into Eq.~(\ref{eq:robustApproxDQ}) we obtain
estimates for the constants $\delta$ and $q$ which in turn will be used
in Eq.~(\ref{eq:robustApprox}) to obtain the approximate distribution
of the robust Wilks' Lambda statistic $\mathit \Lambda_R^{\boldsymbol \cdot}$.

The above procedure to approximate the null distribution of the
robust test statistic depends on the dimension $p$, number of
levels $r$, $c$, and sample size $n$ in each factor combination group
as well as on the considered MANOVA model and the tested hypothesis.
The resulting two parameters of the approximation, $\delta$ and $q$,
can be stored and reused whenever a new MANOVA problem with exactly
the same parameters occurs. Alternatively, we may directly perform
a simulation for the problem at hand in order to compute $p$-values
of the robust tests.

\citet{TODOROV2010} have shown the accuracy of this approximation. 

\section{Monte Carlo simulations}
\label{sec:simulations}

In this section a Monte Carlo study is conducted to assess the
performance of the proposed test 
statistics---the achieved
significance level and the power of the tests. Additionally, 
the behavior of the robust test statistics is evaluated in the presence of
outliers. 

We compare the results of the MCD-based tests to the classical
Wilks' Lambda tests and to an alternative approach proposed by
\citet{NATH1985} which uses the ranks of the observations. 
Although this approach was only proposed for one-way MANOVA models we extend
it to the two-way MANOVA cases. 
The Wilks'
Lambda statistics are calculated on the ranks of the original data and
are referred to as rank transformed statistics,
$\mathit \Lambda_{\mbox{\scriptsize rnk}}^{\boldsymbol \cdot}$.
The distributions of the test statistics are approximated by those of the
normal testing theory.

We considered several dimensions $p \in \{2, 6\}$, number of levels
$r, c \in \{2, 3, 5\}$, and number of samples $n \in \{20, 30, 50\}$
assuming an equal sample size in each factor combination group. The different
designs are given in Table~\ref{tab:designs}. These 18 designs were 
used to simulate data according to both considered models, the two-way MANOVA
model with main effects only as well as the one with main effects and 
interactions. 

\begin{table}[thp]
  \centering
  \caption{Selected designs for the simulation study}
  \label{tab:designs}
  \begin{tabular}{lll rrr}
    \hline
    $r$ & $c$ & $p$ & \multicolumn{3}{c}{$n$} \\
    \hline
    2 & 2 & 2 & 20 & 30 & 50 \\
      &   & 6 & 20 & 30 & 50 \\
    3 & 2 & 2 & 20 & 30 & 50 \\
      &   & 6 & 20 & 30 & 50 \\
    5 & 2 & 2 & 20 & 30 & 50 \\
      &   & 6 & 20 & 30 & 50 \\
    \hline
  \end{tabular}
\end{table}

Here, the significance level $\alpha$ is set to 0.05.

\subsection{Finite-sample accuracy}

Under the null hypotheses 
of $\boldsymbol \alpha_i$,
$\boldsymbol \beta_j$, and $\boldsymbol \gamma_{ij}$ being zero in
both models 
we assume that the
observations come from identical multivariate distributions,
i.e., $\boldsymbol \mu_{11} = \ldots = \boldsymbol \mu_{ij} =
\ldots = \boldsymbol \mu_{rc} = \boldsymbol \mu$, with $i = 1, \ldots, r$
and $j = 1, \ldots, c$.
Since the considered statistics are affine equivariant, without
loss of generality, we can assume $\boldsymbol \mu$ to be a null
vector, i.e., $\boldsymbol \mu = (0, \ldots, 0)^\top = \boldsymbol 0$,  and
the covariance matrix to be $\boldsymbol I_p$. Thus,
for each design listed in Table~\ref{tab:designs} 
we generate
$N = rcn$ $p$-variate vectors distributed as
$\mathcal N_p ( \boldsymbol 0, \boldsymbol I_p )$ and calculate the classical
Wilks' Lambda test statistics, 
$\mathit \Lambda_{\mbox{\scriptsize cla}}^{\boldsymbol \cdot}$, 
the rank transformed ones, 
$\mathit \Lambda_{\mbox{\scriptsize rnk}}^{\boldsymbol \cdot}$, and 
the robust versions, 
$\mathit \Lambda_{R}^{\boldsymbol \cdot}$, based on MCD
estimates.
This is repeated $m = 1000$ times and the percentage
of values of the test statistics above the appropriate critical value
of the corresponding approximate distribution is taken as an estimate
of the true significance level.
The classical Wilks' Lambda and the
rank transformed Wilks' Lambda are compared to the Bartlett
approximation given by Eq.~(\ref{eq:BartlettApprox})  
while the MCD Wilks' Lambda is compared
to the approximate distribution given in Eq.~(\ref{eq:robustApprox})
with parameters $\delta$ and $q$ estimated by Eq.~(\ref{eq:robustApproxDQ}).

\subsection{Finite-sample power comparisons} 

In order to assess the power of the robust Wilks' Lambda statistic we
will generate data samples under an alternative hypothesis
and will examine the frequency of incorrectly failing to reject 
the null hypothesis, i.e., the frequency of Type II errors.
The same combinations of dimensions $p$, number of levels $r$, $c$,
and sample sizes $n$ in each factor combination group as in the experiments
for studying the finite-sample accuracy will be used. There are
infinitely many possibilities for selecting an alternative 
but for the purpose of the study we will borrow an idea from experimental
design: 
we will set the values of the levels according to the {\em least favorable
case\/},
i.e., we will choose the values of the levels in a way for which the null
hypothesis is hardest to reject
resulting in the lowest possible power \citep[cf., e.g., ][]{RASCH2012}.

Moreover, we will distinguish between both considered two-way MANOVA
models. For the two-way MANOVA with interactions the data samples are
generated from a $p$-dimensional normal distribution, where each
factor combination group 
$\mathcal Y_{ij}$, $i = 1, \ldots, r$, $j = 1, \ldots, c$, has a
different mean $\boldsymbol \mu_{ij}$ and all of them have the same
covariance matrix $\boldsymbol I_p$,
\begin{eqnarray}
  \label{eq:simulationWithInter}
  \boldsymbol y_{ijk}
  \sim \mathcal N_p (\boldsymbol \mu_{ij}, \boldsymbol I_p) \ ,
  \quad i = 1, \ldots, r, \quad j = 1, \ldots, c, \quad k = 1, \ldots, n,
\end{eqnarray}
with
\begin{eqnarray*}
  \boldsymbol\mu_{11} = (d/4, 0, \ldots, 0)^\top & , & \quad  
  \boldsymbol\mu_{r1} = (-d/4, 0, \ldots, 0)^\top \ ,\\ 
  \boldsymbol\mu_{1c} = (-d/4, 0, \ldots, 0)^\top & , & \quad  
  \boldsymbol\mu_{rc} = (d/4, 0, \ldots, 0)^\top \ ,\\
  \boldsymbol\mu_{ij} = (0, 0, \ldots, 0)^\top & , & \quad 
  \mbox{for all other $i$'s and $j$'s.}
\end{eqnarray*}
The parameter $d$ takes the following values: $d = 0.0$, 0.2, 0.5, 0.7, 1.0,
1.5, 2.0.
We note that according to the chosen data generating procedure the
null hypotheses $H_A$ and $H_B$ are true for the simulated data. 

Similar, for the two-way MANOVA without interactions the data samples are
generated again from a $p$-dimensional normal distribution, where each
factor combination group 
$\mathcal Y_{ij}$, $i = 1, \ldots, r$, $j = 1, \ldots, c$, has a
different mean $\boldsymbol \mu_{ij}$ and all of them have the same
covariance matrix $\boldsymbol I_p$,
\begin{eqnarray}
  \label{eq:simulationWithoutInter}
  \boldsymbol y_{ijk}
  \sim \mathcal N_p (\boldsymbol \mu_{ij}, \boldsymbol I_p) \ ,
  \quad i = 1, \ldots, r, \quad j = 1, \ldots, c, \quad k = 1, \ldots, n, 
\end{eqnarray}
with
\begin{eqnarray*}
  \boldsymbol\mu_{1j} = (d/2, 0, \ldots, 0)^\top & , & \  
  \boldsymbol\mu_{2j} = (-d/2, 0, \ldots, 0)^\top \ ,\\
  \boldsymbol\mu_{ij} = (0, 0, \ldots, 0)^\top & , & \quad
  i = 3, \ldots, r, \quad j = 1, \ldots, c \ .
\end{eqnarray*}
The parameter $d$ again takes the following values: $d = 0.0$, 0.2, 0.5, 0.7,
1.0, 1.5, 2.0.
We note that according to the chosen data generating procedure the
null hypothesis $H_B$ is true for the simulated data. 

Then we calculate the classical
Wilks' Lambda test statistics, 
$\mathit \Lambda_{\mbox{\scriptsize cla}}^{\boldsymbol \cdot}$, 
the rank transformed ones, 
$\mathit \Lambda_{\mbox{\scriptsize rnk}}^{\boldsymbol \cdot}$, and 
the robust versions, 
$\mathit \Lambda_{R}^{\boldsymbol \cdot}$, based on MCD
estimates.
This is repeated $m = 1000$ times and the rejection frequency where
the statistic exceeds its appropriate critical value is the estimate
of the power for the specific configuration.

\subsection{Finite-sample robustness comparisons} 

Here, for each design in Table~\ref{tab:designs} 
data samples are generated under the null hypotheses for all factor
combination groups $\mathcal Y_{ij}$, $i = 1, \ldots, r - 1$, $j = 1, \ldots, c - 1$, 
i.e., all $\mathcal Y_{ij}$, $i = 1, \ldots, r - 1$, $j = 1, \ldots, c - 1$,  
are distributed as  
$\mathcal N_p (\boldsymbol 0, \boldsymbol I_p)$.
The factor combination group $\mathcal Y_{rc}$
follows the contamination model 
\begin{eqnarray}
  (1 - \varepsilon) \mathcal N_p (\boldsymbol 0, \boldsymbol I_p) + 
  \varepsilon \mathcal N_p (\boldsymbol \mu^*, 0.25^2\boldsymbol I_p) \ ,
\end{eqnarray}
with $\boldsymbol \mu^* = (\nu Q_p, \ldots, \nu Q_p)^\top$ and 
      $Q_p = \sqrt{\chi^2_{p;0.999} / p}$.
The amount of contamination $\varepsilon$ is set to 0.1 and the
outlier distance $\nu$ takes values 2.0, 5.0, 10.0.
By adding $\nu Q_p$ to each component of the outliers we guarantee a
comparable shift for different dimensions \citep[cf.][]{ROCKE1996}.

We then calculate the classical
Wilks' Lambda test statistics, 
$\mathit \Lambda_{\mbox{\scriptsize cla}}^{\boldsymbol \cdot}$, 
the rank transformed ones, 
$\mathit \Lambda_{\mbox{\scriptsize rnk}}^{\boldsymbol \cdot}$, and 
the robust versions, 
$\mathit \Lambda_{R}^{\boldsymbol \cdot}$, based on MCD
estimates.
This is repeated $m = 1000$ times and the percentage
of values of the test statistics above the appropriate critical value
of the corresponding approximate distribution is taken as an estimate
of the true significance level.

\section{Simulation results}
\label{sec:results}

In this section selected results of the simulation study are presented.
Here we only consider two-way MANOVA designs with and without 
interactions with
$r=3$, $c=2$, $n=30$, and $p=2$ for the classical Wilks' Lambda test
(solid line), 
its rank-transformed version (dotted line),
and the MCD-based test (dash-dotted line).
The corresponding results for other designs and dimensions $p$ are similar.
Additional numerical results may be found in \ref{sec:appendix}.
All figures related to the two-way MANOVA with interactions contain
three plots that correspond to the
three possible hypotheses, $H_A$ (left), $H_B$ (middle), and $H_{AB}$
(right), that are able to be tested. 
Whereas all figures related to the two-way MANOVA without interactions
contain two plots that correspond to the
two possible hypotheses, $H_A$ (left) and $H_B$ (right),
that may be tested. 

\subsection{Finite-sample accuracy and power comparisons} 

\begin{figure}[thp]
  \centering
  \includegraphics[width=\textwidth]{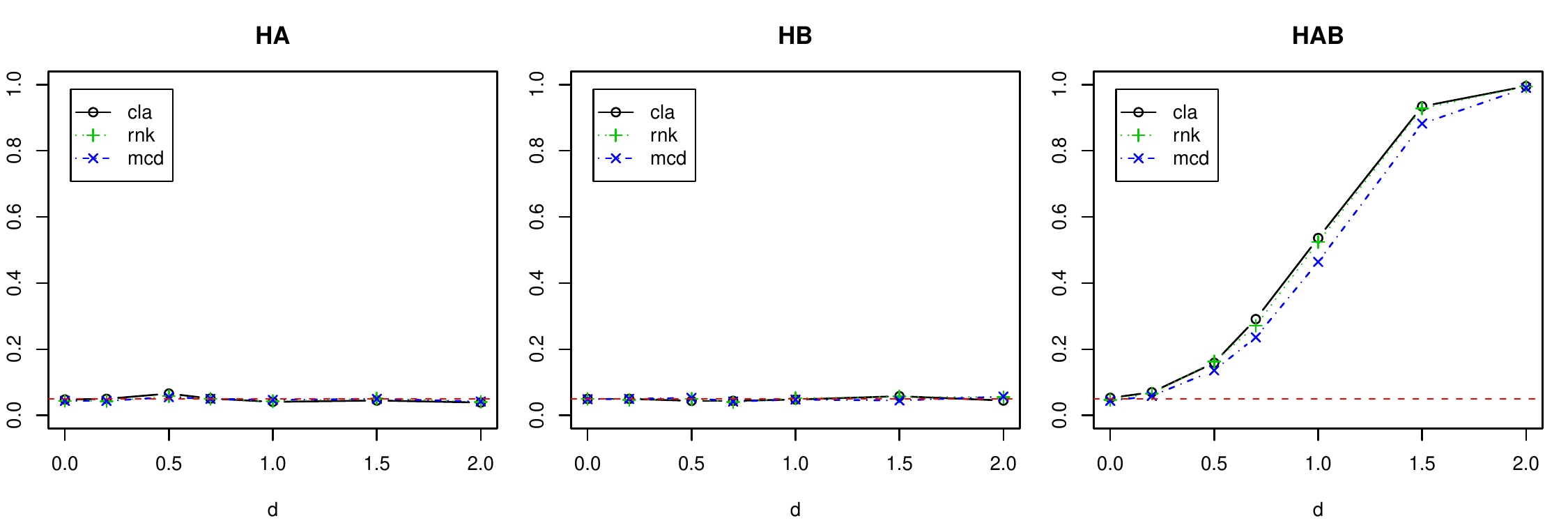}
  \caption{Power of two-way MANOVA with interactions}
  \label{fig:pwrWithInter}
\end{figure}

First, we consider the two-way MANOVA with interactions.
The estimated power is the observed percentage of samples for
which the calculated $p$-values are below 0.05.
Fig.~\ref{fig:pwrWithInter} shows the estimated power
of the two-way MANOVA with interactions. 
Its power curves 
are plotted as a function
of $d$ and include the case $d=0.0$.
The horizontal dashed line indicates the 5\% significance level. 
We note that according to 
Section~\ref{sec:simulations} the null hypotheses $H_A$ and $H_B$ are
true for the simulated data. Hence, in the plots on the left and in the
middle of Fig.~\ref{fig:pwrWithInter}
it is clearly visible that these tests 
keep the significance level independent of the value of $d$.
Further, it can be seen 
in the plot on the right of Fig.~\ref{fig:pwrWithInter}
that the power of the robust test procedure
is almost as high as that of the classical Wilks' Lambda test and its
rank-transformed version.

\begin{figure}[thp]
  \centering
  \includegraphics[width=\textwidth]{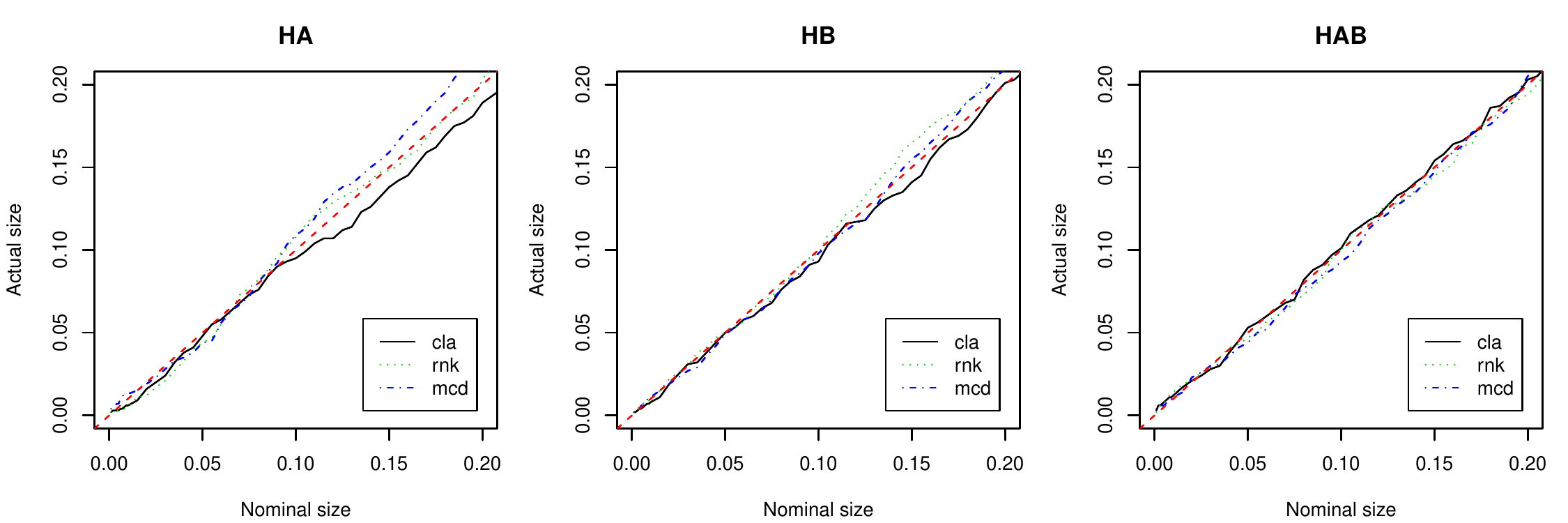}
  \caption{P value plot of two-way MANOVA with interactions}
  \label{fig:pValueWithInter}
\end{figure}

To give a more complete picture of how test statistics follow the
approximate distribution under the null hypothesis in simulated
samples, we will make use of the $P$ value plots proposed by
\citet{DAVIDSON1998}.
Fig.~\ref{fig:pValueWithInter} shows plots of
the empirical distribution functions of the $p$-values of the two-way
MANOVA with interactions. 
The most interesting part of the $P$ value plot is the region where
the size ranges from zero to 0.2 since in practice a significance
level above 20\% is never used. Therefore we limited both axis to
$p$-values $\le 0.2$. We expect that the results in the $P$ value plot
follow the $45^\circ$ line since the $p$-values are distributed
uniformly on $(0, 1)$ if the distribution of the test statistic is
correct.
As it can be seen in Fig.~\ref{fig:pValueWithInter} all test follow
closely the $45^\circ$ line.

\begin{figure}[thp]
  \centering
  \includegraphics[width=\textwidth]{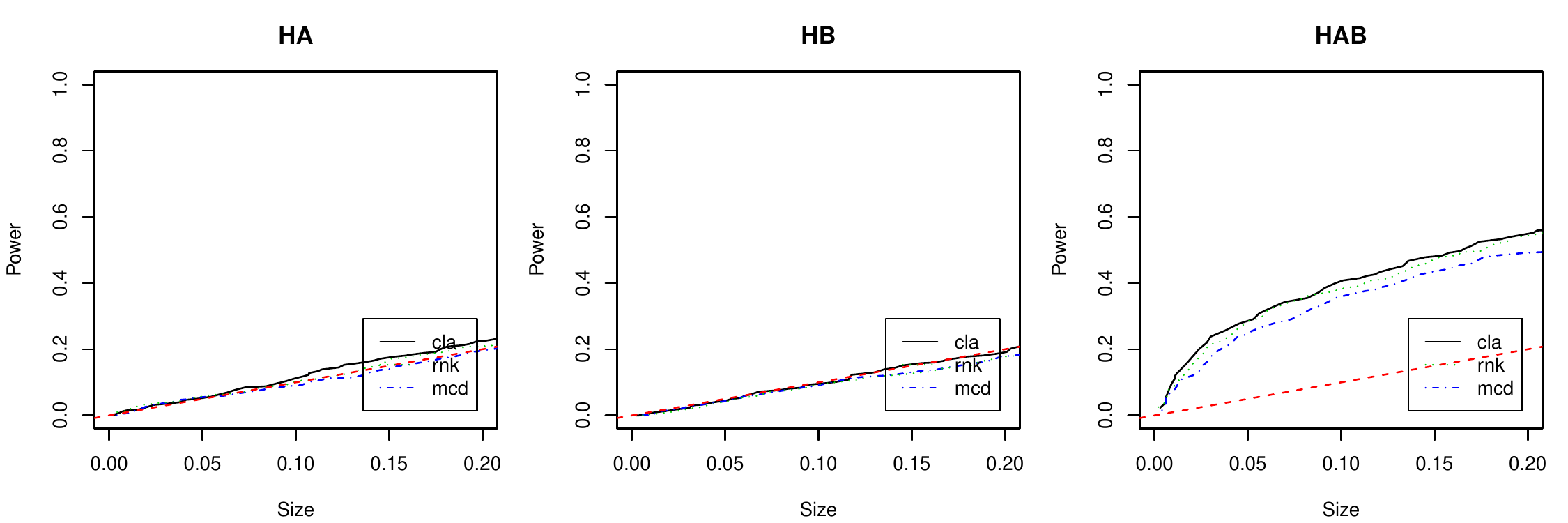}
  \caption{Size-power plot of two-way MANOVA with interactions ($d = 1.0$)}
  \label{fig:sizePwrWithInter}
\end{figure}

Furthermore, the power of different test statistics can be visually
compared by computing size-power curves under fixed alternatives, as
proposed by \citet{DAVIDSON1998}.  As the construction of size-power
plots does not require knowledge of the exact distribution of the test
statistic we may use size-power plots to examine to what extent the
test statistic can differentiate between the null hypothesis and the
alternative.  Fig.~\ref{fig:sizePwrWithInter} shows size-power curves
for three different tests
setting $d = 1.0$. 
The size-power curve should lie above the $45^\circ$ line, the larger the
distance between the curve and the $45^\circ$ line the better.
Again, the most
interesting part of the size-power curve is the region where the size
ranges from zero to 0.2 since in practice a significance level above
20\% is never used.
We again note that here, according to 
Section~\ref{sec:simulations}, the null hypotheses $H_A$ and $H_B$ are
true for the simulated data. Hence, in the plots on the left and in the
middle of Fig.~\ref{fig:sizePwrWithInter}
it can be seen that all curves follow the $45^\circ$ line. 
Further, it is clearly visible in the plot on the right of
Fig.~\ref{fig:sizePwrWithInter} that all curves are far above the
$45^\circ$ line. The classical Wilks' Lambda, being the likelihood
ratio statistic under normality, obviously has the highest power, 
closely followed by its rank-transformed version.
However, the robust MCD-based statistic is almost equally powerful. 

\begin{figure}[thp]
  \centering
  \includegraphics[width=0.7\textwidth]{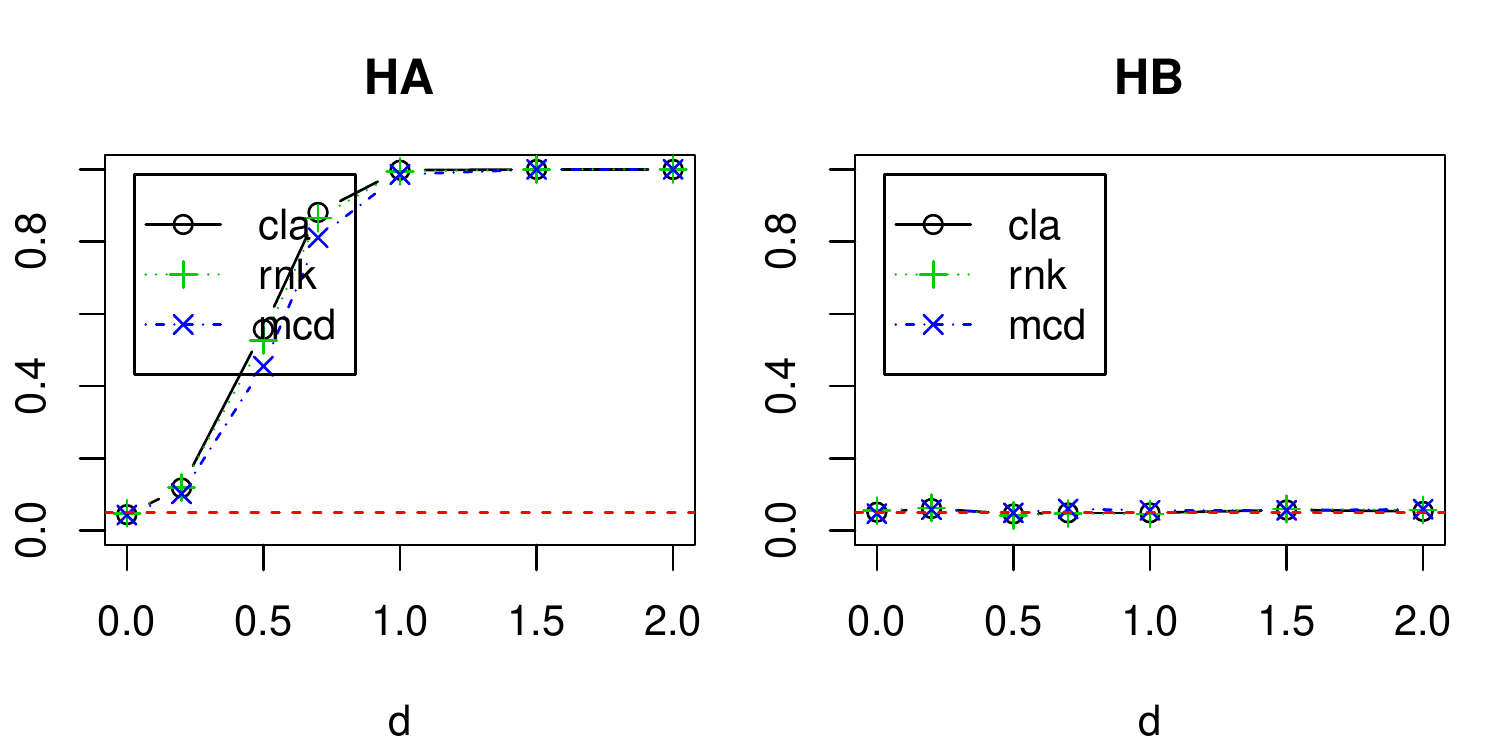}
  \caption{Power of two-way MANOVA without interactions}
  \label{fig:pwrWithoutInter}
\end{figure}

Now, we consider the two-way MANOVA without interactions.
Fig.~\ref{fig:pwrWithoutInter} shows the estimated power
of the two-way MANOVA without interactions. 
Its power curves 
are plotted as a function
of $d$ and include the case $d=0.0$.
The horizontal dashed line again indicates the 5\% significance level. 
We note that according to 
Section~\ref{sec:simulations} the null hypothesis $H_B$ is
true for the simulated data. Hence, in the plot on the right of
Fig.~\ref{fig:pwrWithoutInter}
it is clearly visible that these tests 
keep the significance level independent of the value of $d$.
Further, it can be seen 
in the plot on the left of Fig.~\ref{fig:pwrWithoutInter}
that the power of the robust test procedure
is almost as high as that of the classical Wilks' Lambda test and its
rank-transformed version. 

\begin{figure}[thp]
  \centering
  \includegraphics[width=0.7\textwidth]{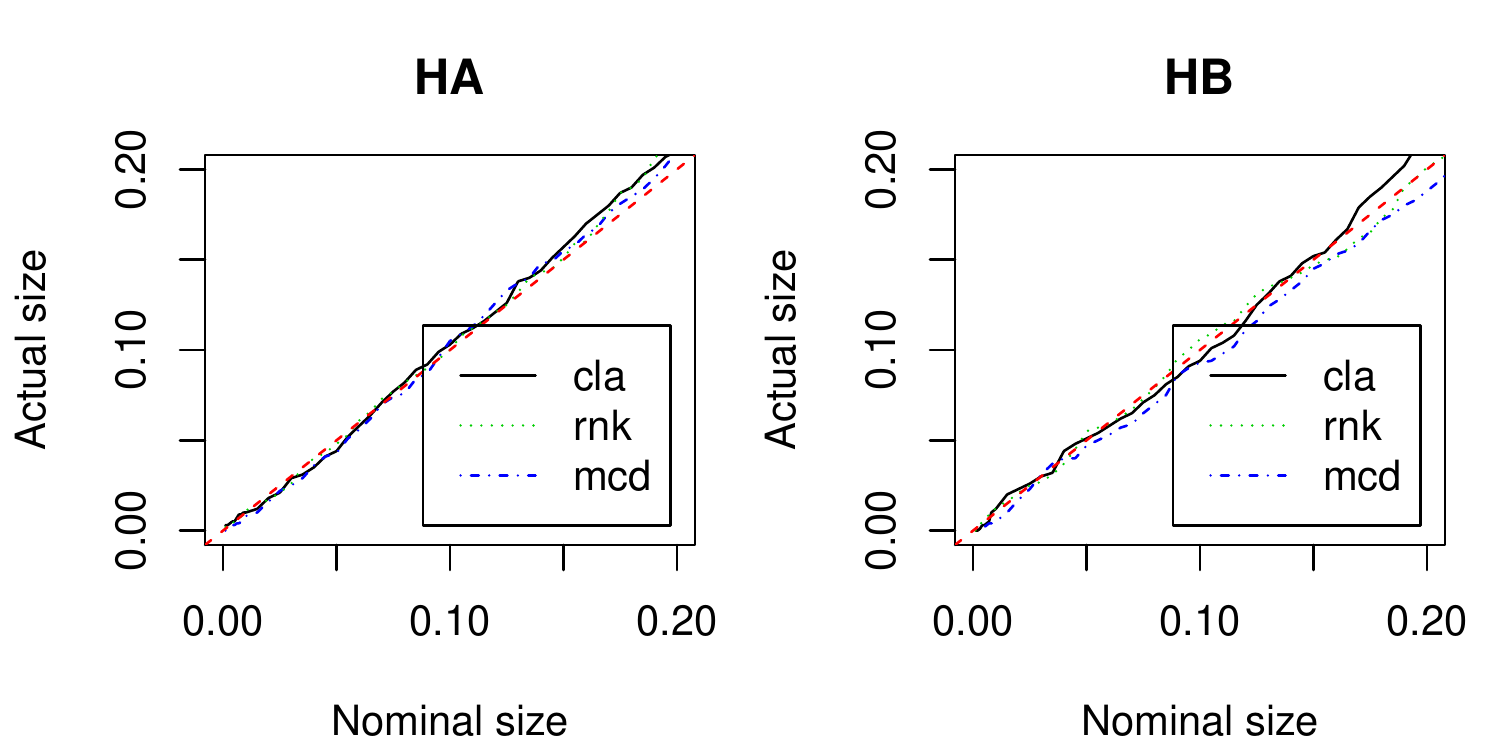}
  \caption{P value plot of two-way MANOVA without interactions}
  \label{fig:pValueWithoutInter}
\end{figure}

Fig.~\ref{fig:pValueWithoutInter} shows $P$ value plots 
of the two-way
MANOVA without interactions. 
As it can be seen in both plots of Fig.~\ref{fig:pValueWithoutInter}
all test follow closely the $45^\circ$ line.

\begin{figure}[thp]
  \centering
  \includegraphics[width=0.7\textwidth]{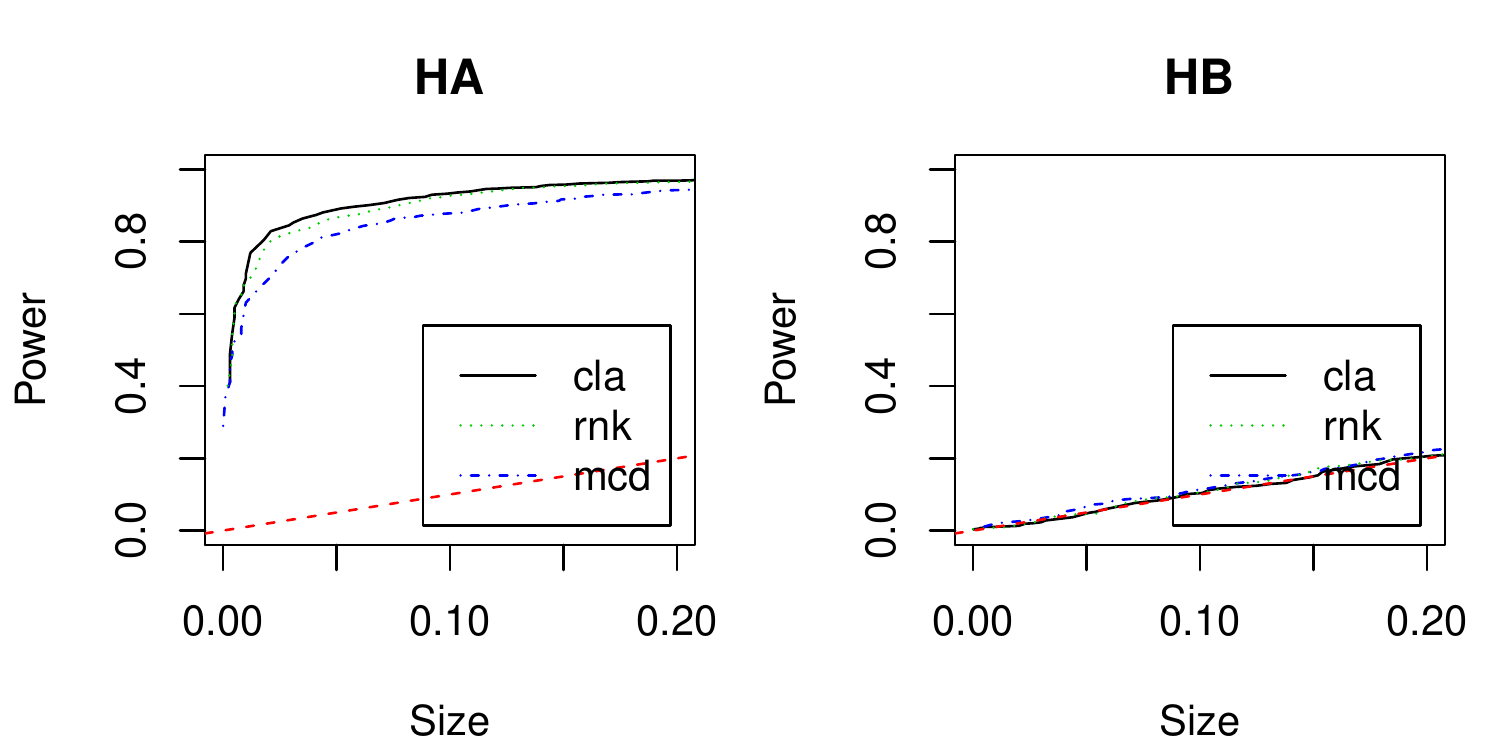}
  \caption{Size-power plot of two-way MANOVA without interactions ($d = 1.0$)}
  \label{fig:sizePwrWithoutInter}
\end{figure}

Fig.~\ref{fig:sizePwrWithoutInter} shows size-power curves
for the three different tests 
setting $d = 1.0$. 
We again note that here, according to 
Section~\ref{sec:simulations}, the null hypothesis $H_B$ is
true for the simulated data. Hence, in the plot on the right of
Fig.~\ref{fig:sizePwrWithoutInter}
it can be seen that all curves follow the $45^\circ$ line. 
Further, it is clearly visible in the plot on the left of
Fig.~\ref{fig:sizePwrWithoutInter} that all curves are far above the
$45^\circ$ line. Again, the classical Wilks' Lambda, being the likelihood
ratio statistic under normality, obviously has the highest power, 
closely followed by its rank-transformed version.
However, the robust MCD-based statistic is almost equally powerful here too. 

\subsection{Finite-sample robustness comparisons}

\begin{figure}[thp]
  \centering
  \includegraphics[width=\textwidth]{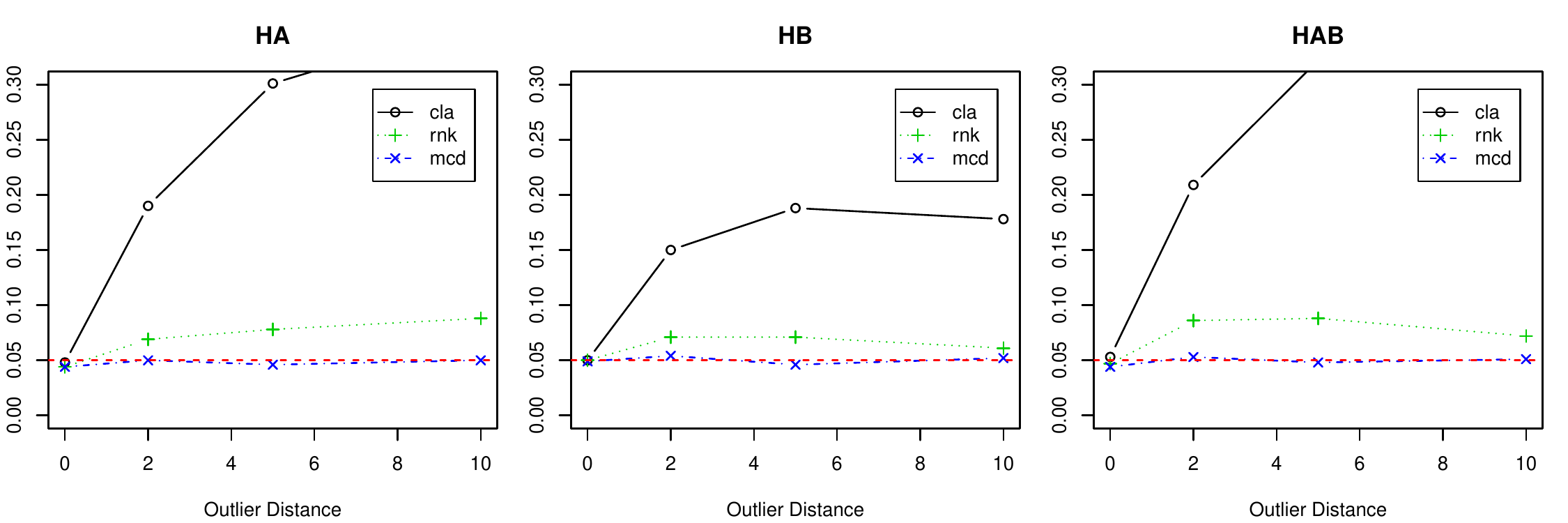}
  \caption{Type I error  of two-way MANOVA with interactions}
  \label{fig:type1WithInter}
\end{figure}

First, we again consider the two-way MANOVA with interactions.
Fig.~\ref{fig:type1WithInter} shows observed Type I error rates of the
two-way MANOVA with interactions 
in the presence of outliers.
The rates are plotted as a function of the outlier distance $\nu$, 
including the non-contaminated case of $\nu = 0$.
The horizontal dashed line indicates the 5\% significance level. 
It is clearly visible in all three plots of Fig.~\ref{fig:type1WithInter}
that the MCD-based test keeps the significance level independent of the
magnitude of the outlier distance $\nu$. Thus, for the MCD-based test,
the Type I error rates turn out to be quite robust and close to the
nominal value.
The rank-transformed test and even more the classical Wilks' Lambda test
are seen to be prone to outliers. 
Both yield very erroneous Type I error rates.

\begin{figure}[tp]
  \centering
  \includegraphics[width=\textwidth]{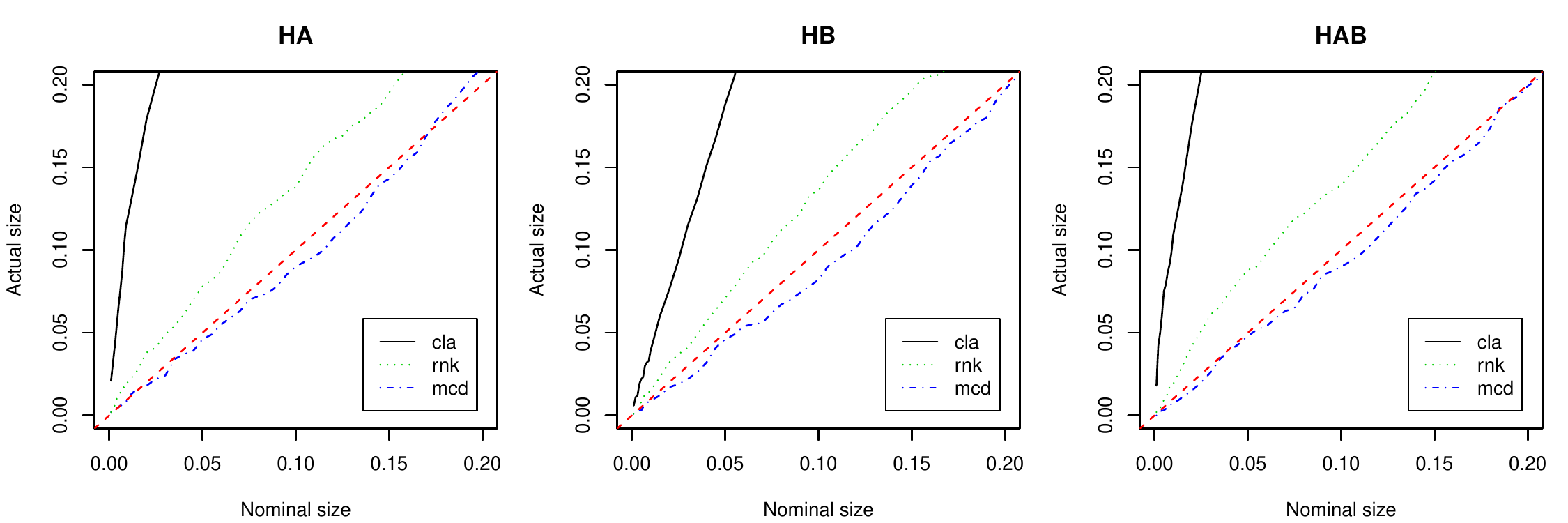}
  \caption{P value plot of two-way MANOVA with interactions (outlier distance 5)}
  \label{fig:pValueWithInterWithOutlier}
\end{figure}

Further, Fig.~\ref{fig:pValueWithInterWithOutlier} shows $P$ value plots
for the three different tests 
in the presence of outliers 
setting the outlier distance $\nu$ to 5.0.
As can be seen in all plots of Fig.~\ref{fig:pValueWithInterWithOutlier}
the MCD-based test follows the $45^\circ$ line fairly well whereas
the classical Wilks' Lambda test and its rank-transformed version
deviates significantly from it. 

\begin{figure}[thp]
  \centering
  \includegraphics[width=0.7\textwidth]{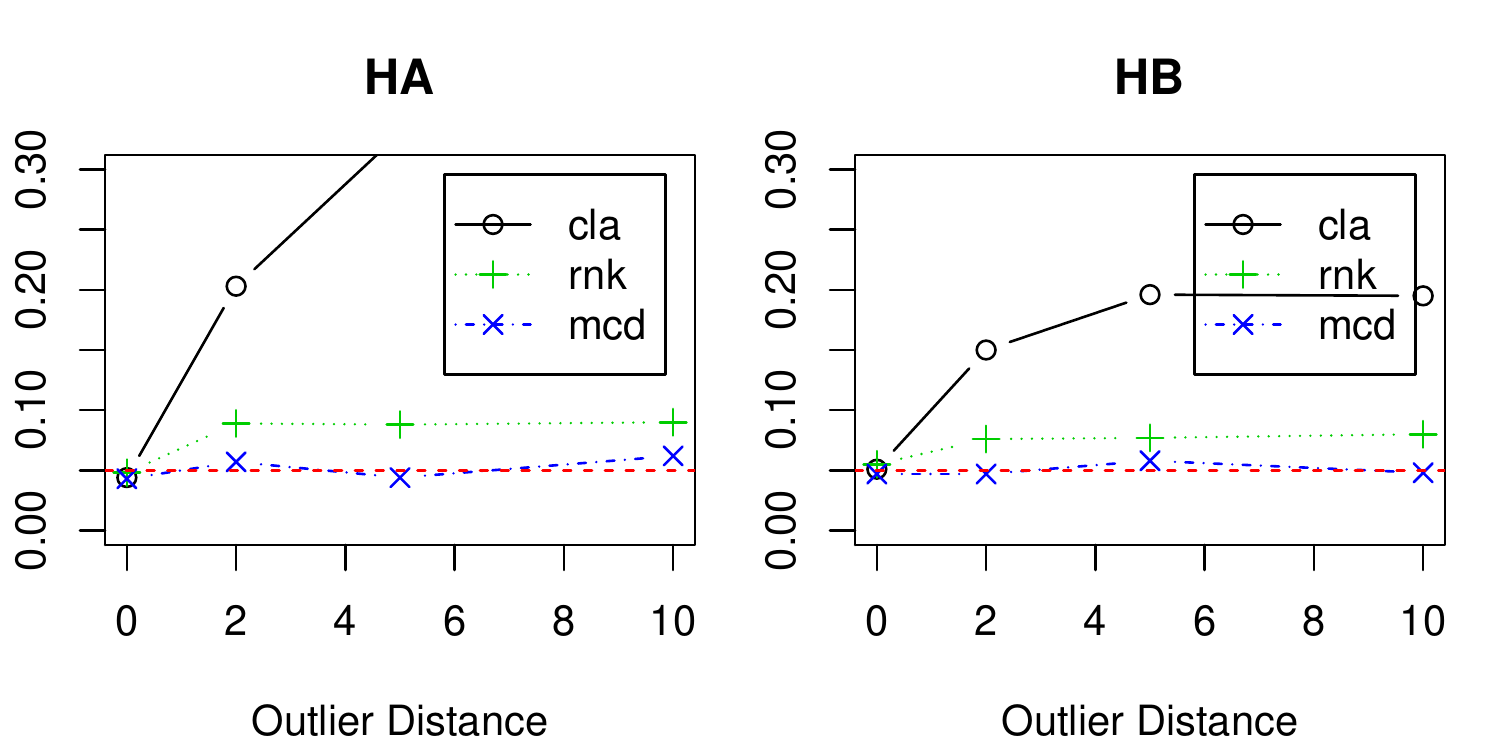}
  \caption{Type I error  of two-way MANOVA without interactions}
  \label{fig:type1WithoutInter}
\end{figure}

Now, we consider the two-way MANOVA without interactions.
Fig.~\ref{fig:type1WithoutInter} shows observed Type I error rates of the
two-way MANOVA without interactions 
in the presence of outliers.
The rates are plotted as a function of the outlier distance $\nu$, 
including the non-contaminated case of $\nu = 0$.
The horizontal dashed line again indicates the 5\% significance level. 
It is clearly visible in both plots of Fig.~\ref{fig:type1WithoutInter}
that the MCD-based test keeps the significance level independent of the
magnitude of the outlier distance $\nu$. Thus again, for the MCD-based test,
the Type I error rates turn out to be quite robust and close to the
nominal value.
The rank-transformed test and even more the classical Wilks' Lambda test
are again seen to be prone to outliers. 
Both yield very erroneous Type I error rates. 

\begin{figure}[thp]
  \centering
  \includegraphics[width=0.7\textwidth]{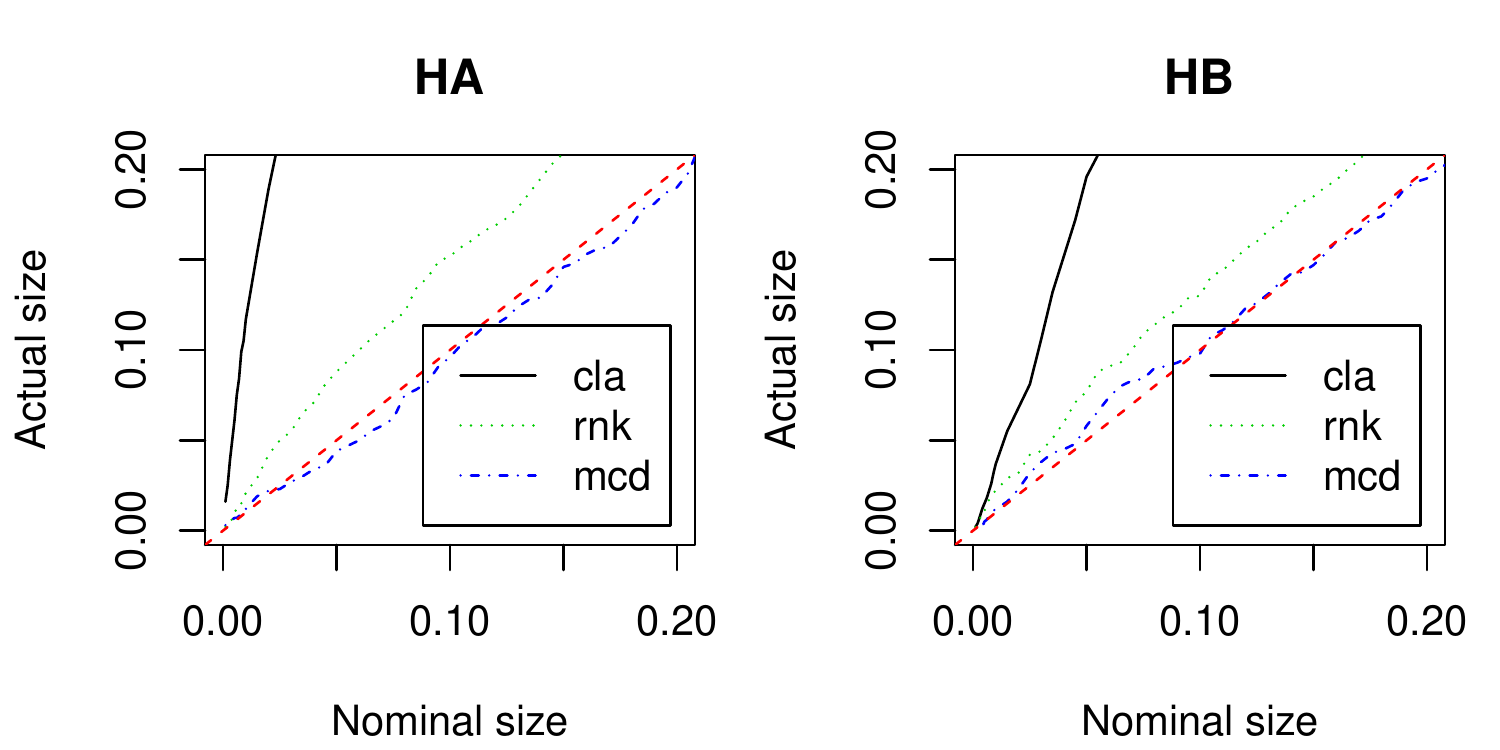}
  \caption{P value plot of two-way MANOVA without interactions
    (outlier distance 5)}
  \label{fig:pValueWithoutInterWithOutlier}
\end{figure}

Further, Fig.~\ref{fig:pValueWithoutInterWithOutlier} shows $P$ value plots
for the three different tests
in the presence of outliers 
setting the outlier distance $\nu$ to 5.0.
As it can be seen in both plots of Fig.~\ref{fig:pValueWithoutInterWithOutlier}
the MCD-based test again follows the $45^\circ$ line fairly well whereas
the classical Wilks' Lambda test and its rank-transformed version
deviates significantly from it. 

\section{Real data example}
\label{sec:example}

We will 
illustrate the application of the proposed robust test statistics
using waste data collected in 
Salzburg, Austria \citep{LEBERSORGER2013}. In the years 2011, 2012, and 2013
waste bins of residual municipal solid waste 
were analyzed in two different districts. The waste
bins were randomly chosen. The numbers of waste bins are given in Table
\ref{tab:numberOfBins}. Due to regulatory reasons the data were anonymized.

\begin{table}[thp]
  \centering
  \caption{Number of analyzed waste bins}
  \label{tab:numberOfBins}
  \begin{tabular}{llrrr}
    \hline
    year     &    & 2011 & 2012 & 2013 \\
    \hline
    district & XY &   30 &   29 &   30 \\
             & A  &   29 &   29 &   29 \\
    \hline                               
  \end{tabular}
\end{table}

\begin{figure}[thp]
  \centering
  \includegraphics[width=\textwidth]{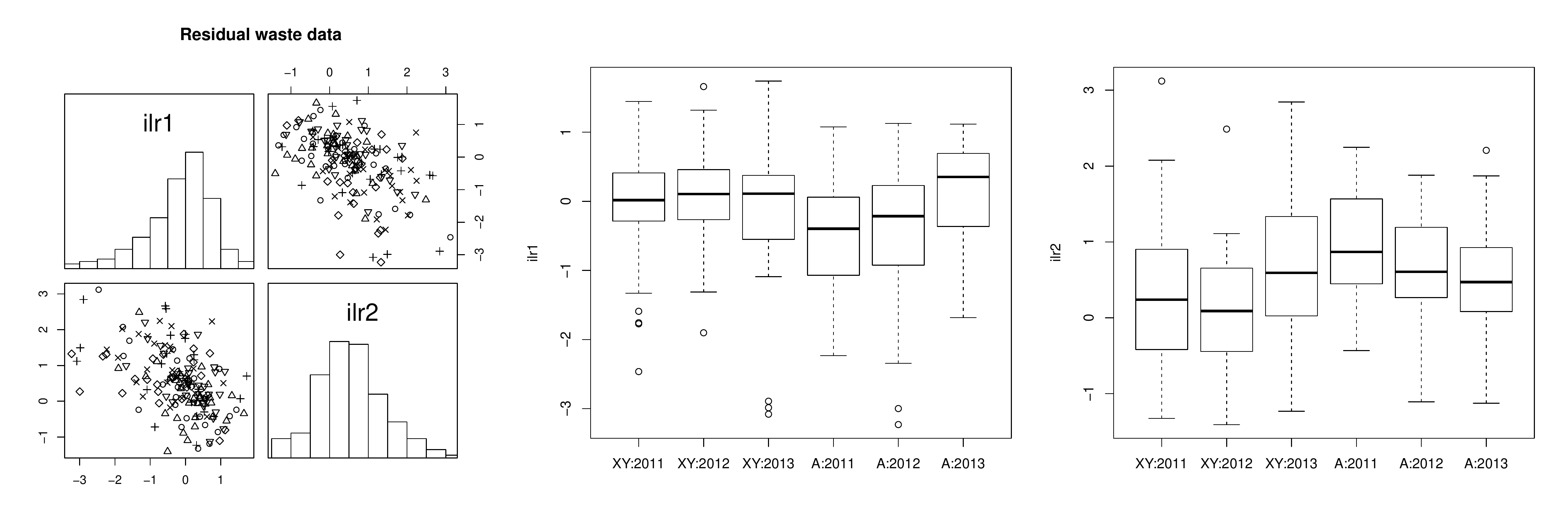}
  \caption{Scatter plot matrix of the isometric log-ratio (ilr)
    transformed waste data (left panel) and grouped boxplots of
    the same data (middle and right panel)}
  \label{fig:groupedBoxplots}
\end{figure}

\begin{figure}[thp]
  \centering
  \includegraphics[width=\textwidth]{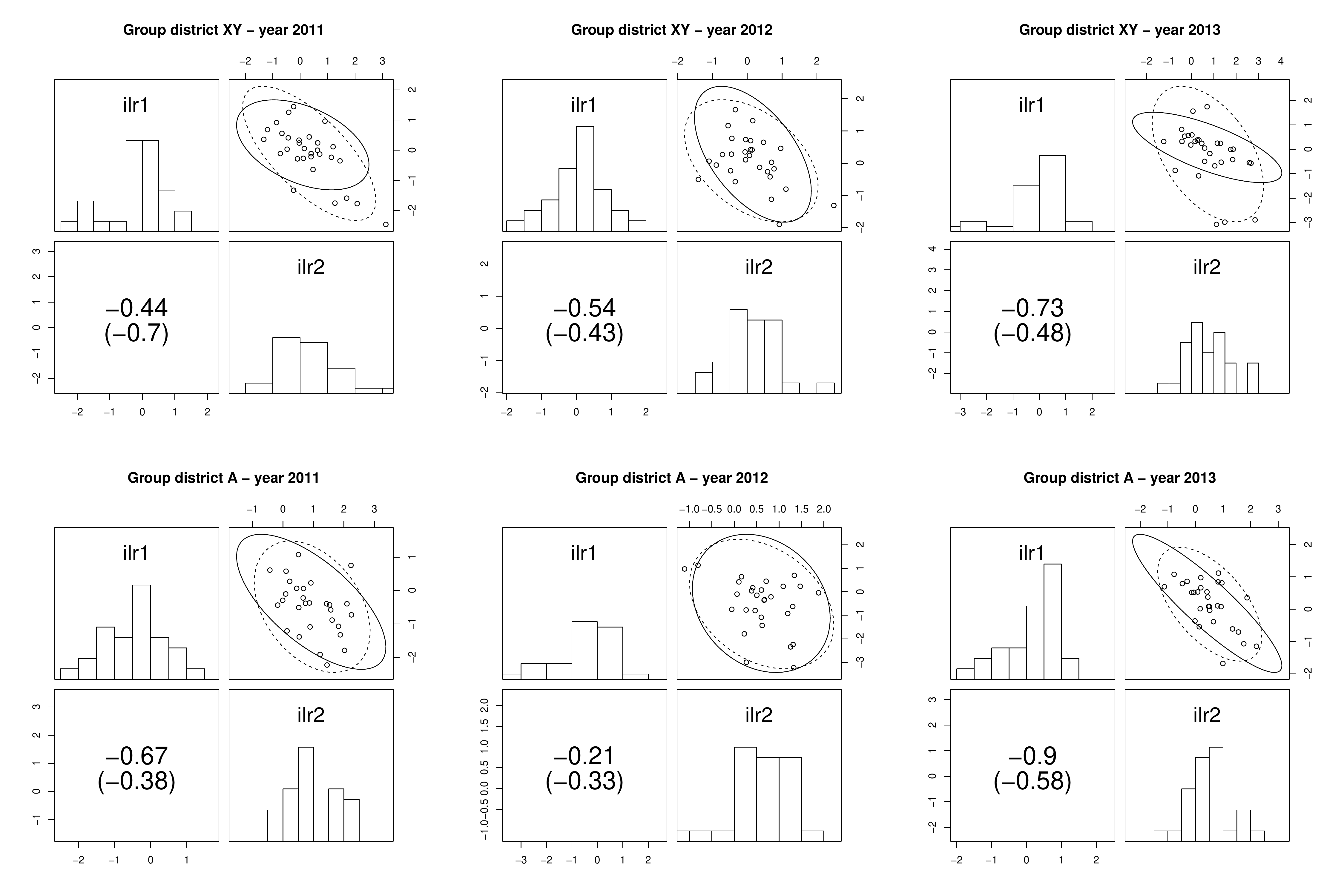}
  \caption{Isometric log-ratio (ilr) transformed waste data: Scatter plot
    matrices of each factor combination group separately with classical
    (dashed lines) and robust (solid lines) 95\% confidence ellipses
    in the upper triangle, classical (in parentheses) and robust
    correlations in the lower triangle and histograms on the diagonal}
  \label{fig:scatterplotMatrices}
\end{figure}

In a residual waste analysis the waste of each bin is divided in up to
20 different fractions and each fraction's portion (given in percentage
of the total weight of the bin) is recorded. For our example we
will only consider 
three main fractions, namely biogenic waste, recyclables, and residual
waste, and aggregate the percentages of corresponding fractions. The first
six observations of the raw data are given in the following:
\begin{verbatim}
  district year biogenic recyclables residual
1       XY 2011   0.2073      0.2493   0.5434
2       XY 2011   0.7065      0.1194   0.1741
3       XY 2011   0.1058      0.6923   0.2019
4       XY 2011   0.2537      0.2985   0.4478
5        A 2011   0.4793      0.1047   0.4160
6        A 2011   0.0966      0.1690   0.7345
...
\end{verbatim}
So, our data matrix consists of $N=176$ rows and $p=3$ columns.
However, we note that as each row sums up to 1 the observations
are compositions being part of the
2-dimensional simplex \citep[cf.][]{AITCHISON1986}.
Most methods from multivariate statistics developed for real valued
data are misleading or inapplicable for compositional data
\citep[cf.][]{BOOGAART2013}.
Hence, we use the isometric log-ratio (ilr) transformation
which is an isometric linear mapping between the $p$-dimensional
simplex and $\mathbb{R}^{p-1}$ 
to obtain a 2-dimensional data matrix for further analysis.
The left panel of Fig.~\ref{fig:groupedBoxplots} shows the scatter plot
matrix of the ilr-transformed data together with histograms of each
variable. The middle and left panels display grouped boxplots for
the different factor combination groups.
In Fig.~\ref{fig:scatterplotMatrices} the scatter plot matrices of each
factor combination group are presented separately with classical
(dashed lines) and robust (solid lines) 95\% confidence ellipses
in the upper triangle, classical (in parentheses) and robust
correlations in the lower triangle and histograms on the diagonal. 

We now perform a two-way MANOVA with and without interactions
using the ilr-transformed waste data.
All possible hypotheses (on row effects, column effects, and, if
applicable, interaction effects) were tested using the classical
Wilks' Lambda test statistic (cla), the rank transformed one (rnk),
and the robust version based on MCD estimates (mcd). 
The $p$-values of the corresponding statistics testing for main
effects only as well as for main effects and interactions are given
in Tables~\ref{tab:pValuesWithout} and \ref{tab:pValuesWith},
respectively.

\begin{table}[htp]
  \centering
  \begin{tabular}{lrrr}
    \hline
    & cla & rnk & mcd \\
    \hline
    district & 0.072 & 0.019 & 0.006 \\
    year     & 0.052 & 0.034 & 0.019 \\
    \hline
  \end{tabular}
  \caption{$p$-values for the classical, rank based, and robust two-way MANOVA tests without interactions applied to the isometric log-ratio (ilr) transformed waste data}
  \label{tab:pValuesWithout}
\end{table}

\begin{table}[htp]
  \centering
  \begin{tabular}{lrrr}
    \hline
    & cla & rnk & mcd \\
    \hline
    district & 0.064 & 0.016 & 0.005 \\
    year     & 0.051 & 0.034 & 0.015 \\
    district:year & 0.019 &  0.028 & 0.054 \\
    \hline
  \end{tabular}
  \caption{$p$-values for the classical, rank based, and robust two-way MANOVA tests with interactions applied to the isometric log-ratio (ilr) transformed waste data}
  \label{tab:pValuesWith}
\end{table}

First, we consider MANOVA tests without interactions.
For the tests based on the classical Wilks' Lambda statistic both
hypotheses testing for main effects cannot be rejected at a significance
level of $\alpha = 0.05$, whereas for the rank and MCD-based tests
we can reject both hypotheses. 
The same is true for MANOVA tests with interactions. However, the
hypothesis testing for interactions is rejected for classical and
rank based test, whereas for the MCD-based test it cannot be rejected.
These results coincide with the practitioners' assumptions 
who expected significant main effects of the factors
\texttt{district} and \texttt{year} but no interaction effect. 

\section{Conclusions}
\label{sec:conclusions}

This paper considered robust test statistics for the two-way MANOVA
model.
Robust versions of the Wilks' Lambda statistics testing main effects
as well as interactions were introduced by replacing the classical
estimates for mean and covariance with robust counterparts
which extends the approach of \citet{TODOROV2010}.
Here, the MCD estimator was chosen and approximate distributions
of the robust test statistics were derived by simulations.
The size and power of the new proposed tests were compared with
the classical and rank transformed Wilks' Lambda tests in Monte
Carlo studies.
Various simulations were performed considering different dimensions,
number of levels, and sample sizes of factor combination groups.
Although only a selection of the results is presented in the paper
these are typical representative outcomes.
Therefore it can be concluded that the significance level of the
robust tests is reasonably precise in case of normally distributed
errors as well as in the presence of outliers.
In the latter case it turns out that the actual size of the robust
tests are in general much closer to the nominal size than the classical
and rank transformed Wilks' Lambda tests.
Furthermore, as indicated by the size-power curves, the robust tests
do not lose much power compared to the classical Wilks' Lambda test.
Further research on extending the work of \citet{VANAELST2011}
to two-way MANOVA designs and comparing the results to the robust
tests introduced here is interesting and warranted.

All computations presented in Section~\ref{sec:results}
as well as the waste data example of Section~\ref{sec:example} 
were performed using the statistical
software environment R \citep{RCORE2018}.
The functions used to
perform the two-way MANOVA tests can be obtained at 
\texttt{http://short.boku.ac.at/km53zz}. 
Moreover, the package \texttt{compositions} \citep{BOOGAART2014} was used
to ilr-transform the waste data.

Here, only test statistics in the context of two-way MANOVA were introduced.
However, by following the same principle of partitioning the total
SSP matrix, these test statistics can easily be extended to many
more complex and higher order MANOVA designs.
This is also a topic for further research.

Moreover, the tests considered here focus on robustness against data
contamination, but still assume that the different groups share a common
covariance structure. Hence, although no improvement of the proposed
tests compared to the classical Wilks' Lambda tests with respect
to robustness against heterogeneity of the covariance structure can
be expected, this deserves further research.
For tests, based on the Wilks' Lambda statistic, the robustness against
heterogeneity of the covariance has been extensively investigated.
\citet{RENCHER1998} gives an overview and concludes that only severe
heterogeneity seriously affects Wilks' Lambda test statistics.
An alternative test statistic, the Pillai's trace statistic, is even
more stable in the presence of heterogeneity of covariances.
In future research robust versions of this test statistic may be
studied and compared to the robust tests introduced here. 

\section*{Acknowledgments}

We are grateful to Sandra Lebersorger and the state government of Salzburg 
for providing the waste data set
as well as to Johannes Tintner and Karl Moder for helpful discussions
and comments on an earlier draft of this paper. 

\appendix

\section{Numerical results of the simulation study}
\label{sec:appendix}

In this section additional numerical results of the simulation study are
presented.
The performance of the MCD-based test (mcd) is compared to the
classical Wilks' Lambda test (cla) and its rank-transformed version (rnk). 
Tables related to the two-way MANOVA with interactions
contain only results of testing for interaction effects
whereas tables related to the two-way MANOVA without interactions
contain only results of testing for row effects.
The significance level $\alpha$ was set to 0.05.
The results for other significance levels were found to be similar. 

\subsection{Finite-sample accuracy and power comparisons} 

First, we consider the two-way MANOVA with interactions.
In Table~\ref{tab:pwrWith} the results of the finite-sample accuracy 
and power comparison 
are given.
In the column entitled $d=0.0$ the observed Type I error
rates for a nominal level of 0.05 
of testing $H_{AB}$ 
are shown.
It is clearly
seen that all tests are capable to keep the significance level for
all investigated designs and dimensions $p$. 
The remaining columns of Table~\ref{tab:pwrWith}
give the estimated power of testing $H_{AB}$
for different values of $d$.
Moreover, we note that the figures printed in bold in Table~\ref{tab:pwrWith}
and the plot on the right in Fig.~\ref{fig:pwrWithInter} correspond to each
other.

\begin{table}[thp]\scriptsize
\centering
\caption{Power of testing $H_{AB}$ of two-way MANOVA with interactions}
\label{tab:pwrWith}
\begin{tabular}{lllll rrrrrrr}
  \hline
  $r$ & $c$ & $p$ & $n$ & & \multicolumn{7}{c}{$d$} \\
  \cline{6-12}
   & & & & & \multicolumn{1}{l}{   0.0} & \multicolumn{1}{l}{   0.2} & \multicolumn{1}{l}{   0.5} & \multicolumn{1}{l}{   0.7} & \multicolumn{1}{l}{   1.0} & \multicolumn{1}{l}{   1.5} & \multicolumn{1}{l}{   2.0} \\ 
   \hline
  2 & 2 & 2 & 20  & cla               & 0.056 & 0.071 & 0.164 & 0.261 & 0.474 & 0.848 & 0.979 \\ 
    &   &   &     & rnk               & 0.054 & 0.076 & 0.158 & 0.244 & 0.450 & 0.847 & 0.979 \\ 
    &   &   &     & mcd               & 0.047 & 0.058 & 0.126 & 0.192 & 0.365 & 0.759 & 0.943 \\ 
  2 & 2 & 6 & 20  & cla               & 0.050 & 0.057 & 0.111 & 0.169 & 0.312 & 0.650 & 0.901 \\ 
    &   &   &     & rnk               & 0.044 & 0.053 & 0.112 & 0.163 & 0.297 & 0.637 & 0.902 \\ 
    &   &   &     & mcd               & 0.059 & 0.054 & 0.079 & 0.119 & 0.225 & 0.470 & 0.780 \\ 
  2 & 2 & 2 & 30  & cla               & 0.059 & 0.057 & 0.199 & 0.378 & 0.693 & 0.959 & 1.000 \\ 
    &   &   &     & rnk               & 0.065 & 0.067 & 0.195 & 0.359 & 0.659 & 0.947 & 1.000 \\ 
    &   &   &     & mcd               & 0.064 & 0.064 & 0.174 & 0.313 & 0.596 & 0.914 & 0.997 \\ 
  2 & 2 & 6 & 30  & cla               & 0.049 & 0.069 & 0.139 & 0.248 & 0.458 & 0.865 & 0.994 \\ 
    &   &   &     & rnk               & 0.049 & 0.068 & 0.147 & 0.243 & 0.457 & 0.847 & 0.990 \\ 
    &   &   &     & mcd               & 0.061 & 0.060 & 0.124 & 0.202 & 0.387 & 0.775 & 0.987 \\ 
  2 & 2 & 2 & 50  & cla               & 0.062 & 0.091 & 0.318 & 0.603 & 0.895 & 0.998 & 1.000 \\ 
    &   &   &     & rnk               & 0.062 & 0.091 & 0.317 & 0.584 & 0.884 & 0.999 & 1.000 \\ 
    &   &   &     & mcd               & 0.059 & 0.086 & 0.284 & 0.533 & 0.844 & 0.996 & 1.000 \\ 
  2 & 2 & 6 & 50  & cla               & 0.055 & 0.062 & 0.212 & 0.405 & 0.753 & 0.993 & 1.000 \\ 
    &   &   &     & rnk               & 0.060 & 0.061 & 0.199 & 0.388 & 0.734 & 0.994 & 1.000 \\ 
    &   &   &     & mcd               & 0.055 & 0.057 & 0.189 & 0.352 & 0.687 & 0.982 & 0.999 \\ 
  3 & 2 & 2 & 20  & cla               & 0.048 & 0.054 & 0.111 & 0.223 & 0.404 & 0.787 & 0.959 \\ 
    &   &   &     & rnk               & 0.051 & 0.055 & 0.117 & 0.205 & 0.394 & 0.765 & 0.949 \\ 
    &   &   &     & mcd               & 0.054 & 0.056 & 0.100 & 0.165 & 0.297 & 0.676 & 0.901 \\ 
  3 & 2 & 6 & 20  & cla               & 0.053 & 0.050 & 0.089 & 0.131 & 0.251 & 0.539 & 0.827 \\ 
    &   &   &     & rnk               & 0.050 & 0.055 & 0.097 & 0.138 & 0.235 & 0.538 & 0.813 \\ 
    &   &   &     & mcd               & 0.043 & 0.052 & 0.074 & 0.109 & 0.185 & 0.424 & 0.706 \\ 
  3 & 2 & 2 & 30  & cla               & \textbf{0.053} & \textbf{0.070} & \textbf{0.159} & \textbf{0.291} & \textbf{0.536} & \textbf{0.934} & \textbf{0.995} \\ 
    &   &   &     & rnk               & \textbf{0.047} & \textbf{0.066} & \textbf{0.163} & \textbf{0.271} & \textbf{0.524} & \textbf{0.928} & \textbf{0.994} \\ 
    &   &   &     & mcd               & \textbf{0.044} & \textbf{0.059} & \textbf{0.136} & \textbf{0.236} & \textbf{0.464} & \textbf{0.882} & \textbf{0.990} \\ 
  3 & 2 & 6 & 30  & cla               & 0.049 & 0.049 & 0.098 & 0.185 & 0.349 & 0.763 & 0.970 \\ 
    &   &   &     & rnk               & 0.048 & 0.044 & 0.099 & 0.180 & 0.333 & 0.746 & 0.956 \\ 
    &   &   &     & mcd               & 0.039 & 0.046 & 0.096 & 0.149 & 0.316 & 0.678 & 0.933 \\ 
  3 & 2 & 2 & 50  & cla               & 0.053 & 0.075 & 0.257 & 0.439 & 0.827 & 0.995 & 1.000 \\ 
    &   &   &     & rnk               & 0.053 & 0.074 & 0.246 & 0.424 & 0.813 & 0.993 & 1.000 \\ 
    &   &   &     & mcd               & 0.058 & 0.069 & 0.219 & 0.379 & 0.752 & 0.986 & 1.000 \\ 
  3 & 2 & 6 & 50  & cla               & 0.058 & 0.056 & 0.133 & 0.272 & 0.596 & 0.953 & 1.000 \\ 
    &   &   &     & rnk               & 0.051 & 0.055 & 0.130 & 0.270 & 0.569 & 0.943 & 1.000 \\ 
    &   &   &     & mcd               & 0.052 & 0.063 & 0.134 & 0.240 & 0.533 & 0.939 & 0.999 \\ 
  5 & 2 & 2 & 20  & cla               & 0.054 & 0.053 & 0.090 & 0.156 & 0.289 & 0.636 & 0.906 \\ 
    &   &   &     & rnk               & 0.047 & 0.049 & 0.098 & 0.149 & 0.282 & 0.613 & 0.893 \\ 
    &   &   &     & mcd               & 0.051 & 0.065 & 0.086 & 0.131 & 0.255 & 0.527 & 0.837 \\ 
  5 & 2 & 6 & 20  & cla               & 0.058 & 0.048 & 0.069 & 0.116 & 0.179 & 0.378 & 0.715 \\ 
    &   &   &     & rnk               & 0.062 & 0.054 & 0.066 & 0.108 & 0.158 & 0.381 & 0.693 \\ 
    &   &   &     & mcd               & 0.046 & 0.066 & 0.068 & 0.096 & 0.148 & 0.325 & 0.592 \\ 
  5 & 2 & 2 & 30  & cla               & 0.055 & 0.059 & 0.123 & 0.240 & 0.442 & 0.847 & 0.987 \\ 
    &   &   &     & rnk               & 0.053 & 0.057 & 0.122 & 0.242 & 0.421 & 0.825 & 0.985 \\ 
    &   &   &     & mcd               & 0.067 & 0.056 & 0.103 & 0.201 & 0.359 & 0.767 & 0.969 \\ 
  5 & 2 & 6 & 30  & cla               & 0.055 & 0.058 & 0.083 & 0.128 & 0.246 & 0.609 & 0.918 \\ 
    &   &   &     & rnk               & 0.051 & 0.063 & 0.082 & 0.126 & 0.243 & 0.591 & 0.888 \\ 
    &   &   &     & mcd               & 0.053 & 0.056 & 0.071 & 0.106 & 0.220 & 0.540 & 0.860 \\ 
  5 & 2 & 2 & 50  & cla               & 0.053 & 0.082 & 0.189 & 0.352 & 0.706 & 0.981 & 0.999 \\ 
    &   &   &     & rnk               & 0.054 & 0.068 & 0.172 & 0.340 & 0.693 & 0.976 & 0.999 \\ 
    &   &   &     & mcd               & 0.059 & 0.057 & 0.181 & 0.299 & 0.639 & 0.958 & 0.999 \\ 
  5 & 2 & 6 & 50  & cla               & 0.044 & 0.070 & 0.134 & 0.227 & 0.462 & 0.889 & 0.997 \\ 
    &   &   &     & rnk               & 0.052 & 0.065 & 0.127 & 0.214 & 0.428 & 0.876 & 0.994 \\ 
    &   &   &     & mcd               & 0.050 & 0.061 & 0.124 & 0.206 & 0.401 & 0.825 & 0.994 \\ 
   \hline
\end{tabular}
\end{table}

Now, we consider the two-way MANOVA without interactions.
In Table~\ref{tab:pwrWithout} the results of the finite-sample accuracy 
and power comparison 
are given as before. 
In the column entitled $d=0.0$ the observed Type I error
rates for a nominal level of 0.05 
of testing $H_{A}$ 
are shown.
It is clearly
seen that all tests are capable to keep the significance level for
all investigated designs and dimensions $p$. 
The remaining columns of Table~\ref{tab:pwrWithout}
give the estimated power of testing $H_{A}$ 
for different values of $d$.
Moreover, we note that the figures printed in bold in Table~\ref{tab:pwrWithout}
and the plot on the left in Fig.~\ref{fig:pwrWithoutInter} correspond to each
other.

\begin{table}[thp]\scriptsize
\centering
\caption{Power of testing $H_{A}$ of two-way MANOVA without interactions}
\label{tab:pwrWithout}
\begin{tabular}{lllll rrrrrrr}
  \hline
  $r$ & $c$ & $p$ & $n$ & & \multicolumn{7}{c}{$d$} \\
  \cline{6-12}
   & & & & & \multicolumn{1}{l}{   0.0} & \multicolumn{1}{l}{   0.2} & \multicolumn{1}{l}{   0.5} & \multicolumn{1}{l}{   0.7} & \multicolumn{1}{l}{   1.0} & \multicolumn{1}{l}{   1.5} & \multicolumn{1}{l}{   2.0} \\ 
   \hline
  2 & 2 & 2 & 20  & cla               & 0.049 & 0.121 & 0.508 & 0.796 & 0.979 & 1.000 & 1.000 \\ 
    &   &   &     & rnk               & 0.047 & 0.111 & 0.487 & 0.790 & 0.973 & 1.000 & 1.000 \\ 
    &   &   &     & mcd               & 0.051 & 0.094 & 0.388 & 0.694 & 0.942 & 1.000 & 1.000 \\ 
  2 & 2 & 6 & 20  & cla               & 0.064 & 0.069 & 0.315 & 0.592 & 0.918 & 0.999 & 1.000 \\ 
    &   &   &     & rnk               & 0.066 & 0.068 & 0.304 & 0.588 & 0.901 & 0.999 & 1.000 \\ 
    &   &   &     & mcd               & 0.054 & 0.055 & 0.203 & 0.407 & 0.765 & 0.993 & 1.000 \\ 
  2 & 2 & 2 & 30  & cla               & 0.053 & 0.165 & 0.669 & 0.930 & 0.999 & 1.000 & 1.000 \\ 
    &   &   &     & rnk               & 0.052 & 0.165 & 0.649 & 0.916 & 0.999 & 1.000 & 1.000 \\ 
    &   &   &     & mcd               & 0.046 & 0.131 & 0.586 & 0.876 & 0.996 & 1.000 & 1.000 \\ 
  2 & 2 & 6 & 30  & cla               & 0.047 & 0.089 & 0.456 & 0.792 & 0.993 & 1.000 & 1.000 \\ 
    &   &   &     & rnk               & 0.045 & 0.083 & 0.444 & 0.778 & 0.991 & 1.000 & 1.000 \\ 
    &   &   &     & mcd               & 0.043 & 0.084 & 0.357 & 0.695 & 0.973 & 1.000 & 1.000 \\ 
  2 & 2 & 2 & 50  & cla               & 0.051 & 0.237 & 0.889 & 0.993 & 1.000 & 1.000 & 1.000 \\ 
    &   &   &     & rnk               & 0.049 & 0.219 & 0.873 & 0.992 & 1.000 & 1.000 & 1.000 \\ 
    &   &   &     & mcd               & 0.049 & 0.193 & 0.833 & 0.990 & 1.000 & 1.000 & 1.000 \\ 
  2 & 2 & 6 & 50  & cla               & 0.058 & 0.139 & 0.745 & 0.963 & 1.000 & 1.000 & 1.000 \\ 
    &   &   &     & rnk               & 0.061 & 0.146 & 0.714 & 0.956 & 1.000 & 1.000 & 1.000 \\ 
    &   &   &     & mcd               & 0.052 & 0.128 & 0.677 & 0.941 & 1.000 & 1.000 & 1.000 \\ 
  3 & 2 & 2 & 20  & cla               & 0.054 & 0.084 & 0.377 & 0.691 & 0.954 & 1.000 & 1.000 \\ 
    &   &   &     & rnk               & 0.053 & 0.086 & 0.356 & 0.660 & 0.944 & 1.000 & 1.000 \\ 
    &   &   &     & mcd               & 0.056 & 0.085 & 0.280 & 0.587 & 0.894 & 1.000 & 1.000 \\ 
  3 & 2 & 6 & 20  & cla               & 0.050 & 0.074 & 0.226 & 0.459 & 0.825 & 0.996 & 1.000 \\ 
    &   &   &     & rnk               & 0.053 & 0.075 & 0.222 & 0.437 & 0.817 & 0.997 & 1.000 \\ 
    &   &   &     & mcd               & 0.055 & 0.072 & 0.169 & 0.334 & 0.700 & 0.985 & 1.000 \\ 
  3 & 2 & 2 & 30  & cla               & \textbf{0.044} & \textbf{0.118} & \textbf{0.557} & \textbf{0.881} & \textbf{0.998} & \textbf{1.000} & \textbf{1.000} \\ 
    &   &   &     & rnk               & \textbf{0.048} & \textbf{0.119} & \textbf{0.527} & \textbf{0.865} & \textbf{0.995} & \textbf{1.000} & \textbf{1.000} \\ 
    &   &   &     & mcd               & \textbf{0.043} & \textbf{0.102} & \textbf{0.455} & \textbf{0.811} & \textbf{0.986} & \textbf{1.000} & \textbf{1.000} \\ 
  3 & 2 & 6 & 30  & cla               & 0.046 & 0.080 & 0.360 & 0.709 & 0.975 & 1.000 & 1.000 \\ 
    &   &   &     & rnk               & 0.051 & 0.087 & 0.336 & 0.681 & 0.968 & 1.000 & 1.000 \\ 
    &   &   &     & mcd               & 0.051 & 0.086 & 0.308 & 0.607 & 0.929 & 1.000 & 1.000 \\ 
  3 & 2 & 2 & 50  & cla               & 0.049 & 0.160 & 0.812 & 0.987 & 1.000 & 1.000 & 1.000 \\ 
    &   &   &     & rnk               & 0.051 & 0.160 & 0.791 & 0.982 & 1.000 & 1.000 & 1.000 \\ 
    &   &   &     & mcd               & 0.051 & 0.149 & 0.740 & 0.971 & 1.000 & 1.000 & 1.000 \\ 
  3 & 2 & 6 & 50  & cla               & 0.044 & 0.100 & 0.627 & 0.926 & 0.998 & 1.000 & 1.000 \\ 
    &   &   &     & rnk               & 0.042 & 0.104 & 0.596 & 0.910 & 0.997 & 1.000 & 1.000 \\ 
    &   &   &     & mcd               & 0.047 & 0.090 & 0.576 & 0.884 & 0.997 & 1.000 & 1.000 \\ 
  5 & 2 & 2 & 20  & cla               & 0.045 & 0.079 & 0.296 & 0.562 & 0.906 & 0.999 & 1.000 \\ 
    &   &   &     & rnk               & 0.045 & 0.074 & 0.272 & 0.554 & 0.886 & 1.000 & 1.000 \\ 
    &   &   &     & mcd               & 0.045 & 0.071 & 0.233 & 0.468 & 0.814 & 0.996 & 1.000 \\ 
  5 & 2 & 6 & 20  & cla               & 0.061 & 0.066 & 0.172 & 0.347 & 0.684 & 0.990 & 1.000 \\ 
    &   &   &     & rnk               & 0.055 & 0.069 & 0.163 & 0.332 & 0.662 & 0.986 & 1.000 \\ 
    &   &   &     & mcd               & 0.057 & 0.063 & 0.143 & 0.291 & 0.560 & 0.972 & 1.000 \\ 
  5 & 2 & 2 & 30  & cla               & 0.037 & 0.083 & 0.412 & 0.758 & 0.990 & 1.000 & 1.000 \\ 
    &   &   &     & rnk               & 0.035 & 0.083 & 0.400 & 0.762 & 0.986 & 1.000 & 1.000 \\ 
    &   &   &     & mcd               & 0.035 & 0.073 & 0.333 & 0.680 & 0.976 & 1.000 & 1.000 \\ 
  5 & 2 & 6 & 30  & cla               & 0.047 & 0.066 & 0.263 & 0.545 & 0.912 & 1.000 & 1.000 \\ 
    &   &   &     & rnk               & 0.051 & 0.067 & 0.245 & 0.518 & 0.899 & 1.000 & 1.000 \\ 
    &   &   &     & mcd               & 0.062 & 0.065 & 0.221 & 0.475 & 0.848 & 1.000 & 1.000 \\ 
  5 & 2 & 2 & 50  & cla               & 0.043 & 0.130 & 0.688 & 0.960 & 1.000 & 1.000 & 1.000 \\ 
    &   &   &     & rnk               & 0.047 & 0.131 & 0.662 & 0.946 & 1.000 & 1.000 & 1.000 \\ 
    &   &   &     & mcd               & 0.039 & 0.118 & 0.610 & 0.926 & 1.000 & 1.000 & 1.000 \\ 
  5 & 2 & 6 & 50  & cla               & 0.049 & 0.074 & 0.472 & 0.833 & 0.996 & 1.000 & 1.000 \\ 
    &   &   &     & rnk               & 0.049 & 0.074 & 0.456 & 0.821 & 0.995 & 1.000 & 1.000 \\ 
    &   &   &     & mcd               & 0.046 & 0.065 & 0.410 & 0.773 & 0.986 & 1.000 & 1.000 \\ 
   \hline
\end{tabular}
\end{table}

\subsection{Finite-sample robustness comparisons}

First, we again consider the two-way MANOVA with interactions.
In Table~\ref{tab:type1With} the results of the finite-sample 
robustness comparison 
are given.
We state the observed Type I error rates for a nominal level
of 0.05 of testing $H_{AB}$ 
in the presence of out outliers. Different outlier distances $\nu$
were considered. It is clearly seen that the robust MCD-based test
is capable to keep the significance level for all investigated
designs and dimensions $p$ whereas the classical Wilks' Lambda test
and its rank-transformed version fails to keep the significance
level.
Moreover, we note that the figures printed in bold in
Table~\ref{tab:type1With} and the plot on the right in
Fig.~\ref{fig:type1WithInter} correspond to each other.

\begin{table}[thp]\scriptsize
\centering
\caption{Type I error of testing $H_{AB}$ of two-way MANOVA with interactions}
\label{tab:type1With}
\begin{tabular}{llll rrrrrrrrr}
  \hline
  $r$ & $c$ & $p$ & $n$ & \multicolumn{1}{l}{   cla} & \multicolumn{1}{l}{      } & \multicolumn{1}{l}{      } & \multicolumn{1}{l}{   rnk} & \multicolumn{1}{l}{      } & \multicolumn{1}{l}{      } & \multicolumn{1}{l}{   mcd} & \multicolumn{1}{l}{      } & \multicolumn{1}{l}{      } \\
  \cline{5-7} \cline{8-10} \cline{11-13}
   & & &            & \multicolumn{1}{l}{   2.0} & \multicolumn{1}{l}{   5.0} & \multicolumn{1}{l}{  10.0} & \multicolumn{1}{l}{   2.0} & \multicolumn{1}{l}{   5.0} & \multicolumn{1}{l}{  10.0} & \multicolumn{1}{l}{   2.0} & \multicolumn{1}{l}{   5.0} & \multicolumn{1}{l}{  10.0} \\ 
   \hline
  2 & 2 & 2 & 20             & 0.116 & 0.104 & 0.105 & 0.078 & 0.061 & 0.078 & 0.054 & 0.046 & 0.056 \\ 
  2 & 2 & 6 & 20             & 0.090 & 0.071 & 0.081 & 0.073 & 0.066 & 0.071 & 0.073 & 0.040 & 0.041 \\ 
  2 & 2 & 2 & 30             & 0.152 & 0.204 & 0.222 & 0.085 & 0.080 & 0.075 & 0.042 & 0.048 & 0.036 \\ 
  2 & 2 & 6 & 30             & 0.127 & 0.126 & 0.141 & 0.099 & 0.097 & 0.092 & 0.059 & 0.033 & 0.057 \\ 
  2 & 2 & 2 & 50             & 0.304 & 0.471 & 0.501 & 0.117 & 0.095 & 0.115 & 0.049 & 0.041 & 0.062 \\ 
  2 & 2 & 6 & 50             & 0.195 & 0.262 & 0.278 & 0.106 & 0.123 & 0.122 & 0.045 & 0.046 & 0.050 \\ 
  3 & 2 & 2 & 20             & 0.136 & 0.184 & 0.190 & 0.062 & 0.065 & 0.073 & 0.045 & 0.050 & 0.036 \\ 
  3 & 2 & 6 & 20             & 0.108 & 0.115 & 0.114 & 0.076 & 0.081 & 0.082 & 0.058 & 0.021 & 0.011 \\ 
  3 & 2 & 2 & 30             & \textbf{0.209} & \textbf{0.322} & \textbf{0.354} & \textbf{0.086} & \textbf{0.088} & \textbf{0.072} & \textbf{0.053} & \textbf{0.048} & \textbf{0.051} \\ 
  3 & 2 & 6 & 30             & 0.140 & 0.212 & 0.228 & 0.080 & 0.089 & 0.093 & 0.062 & 0.025 & 0.065 \\ 
  3 & 2 & 2 & 50             & 0.358 & 0.628 & 0.667 & 0.115 & 0.121 & 0.102 & 0.055 & 0.053 & 0.051 \\ 
  3 & 2 & 6 & 50             & 0.274 & 0.435 & 0.438 & 0.118 & 0.138 & 0.133 & 0.065 & 0.063 & 0.057 \\ 
  5 & 2 & 2 & 20             & 0.139 & 0.259 & 0.311 & 0.057 & 0.064 & 0.067 & 0.034 & 0.048 & 0.047 \\ 
  5 & 2 & 6 & 20             & 0.130 & 0.174 & 0.183 & 0.084 & 0.074 & 0.084 & 0.053 & 0.029 & 0.025 \\ 
  5 & 2 & 2 & 30             & 0.234 & 0.428 & 0.510 & 0.075 & 0.082 & 0.087 & 0.047 & 0.055 & 0.044 \\ 
  5 & 2 & 6 & 30             & 0.170 & 0.308 & 0.333 & 0.090 & 0.090 & 0.084 & 0.053 & 0.029 & 0.055 \\ 
  5 & 2 & 2 & 50             & 0.420 & 0.719 & 0.790 & 0.123 & 0.109 & 0.106 & 0.055 & 0.046 & 0.054 \\ 
  5 & 2 & 6 & 50             & 0.342 & 0.566 & 0.638 & 0.122 & 0.124 & 0.127 & 0.064 & 0.042 & 0.050 \\ 
   \hline
\end{tabular}
\end{table}

Now, we consider the two-way MANOVA without interactions.
In Table~\ref{tab:type1Without} the results of the finite-sample 
robustness comparison 
are given.
We state here the observed Type I error rates for a nominal level
of 0.05 of testing $H_{A}$ 
in the presence of out outliers. Again, different outlier distances $\nu$
were considered. It is clearly seen that the robust MCD-based test
is capable to keep the significance level for all investigated
designs and dimensions $p$ whereas the classical Wilks' Lambda test
and its rank-transformed version fails to keep the significance
level.
Moreover, we note that the figures printed in bold in
Table~\ref{tab:type1Without} and the plot on the left in
Fig.~\ref{fig:type1WithoutInter} correspond to each other.

\begin{table}[thp]\scriptsize
\centering
\caption{Type I error of testing $H_{A}$ of two-way MANOVA without interactions}
\label{tab:type1Without}
\begin{tabular}{llll rrrrrrrrr}
  \hline
  $r$ & $c$ & $p$ & $n$ & \multicolumn{1}{l}{   cla} & \multicolumn{1}{l}{      } & \multicolumn{1}{l}{      } & \multicolumn{1}{l}{   rnk} & \multicolumn{1}{l}{      } & \multicolumn{1}{l}{      } & \multicolumn{1}{l}{   mcd} & \multicolumn{1}{l}{      } & \multicolumn{1}{l}{      } \\
  \cline{5-7} \cline{8-10} \cline{11-13}
   & & &            & \multicolumn{1}{l}{   2.0} & \multicolumn{1}{l}{   5.0} & \multicolumn{1}{l}{  10.0} & \multicolumn{1}{l}{   2.0} & \multicolumn{1}{l}{   5.0} & \multicolumn{1}{l}{  10.0} & \multicolumn{1}{l}{   2.0} & \multicolumn{1}{l}{   5.0} & \multicolumn{1}{l}{  10.0} \\ 
   \hline
  2 & 2 & 2 & 20             & 0.102 & 0.097 & 0.090 & 0.067 & 0.064 & 0.081 & 0.038 & 0.044 & 0.063 \\ 
  2 & 2 & 6 & 20             & 0.088 & 0.075 & 0.080 & 0.076 & 0.067 & 0.071 & 0.063 & 0.050 & 0.066 \\ 
  2 & 2 & 2 & 30             & 0.153 & 0.216 & 0.202 & 0.078 & 0.087 & 0.094 & 0.051 & 0.053 & 0.047 \\ 
  2 & 2 & 6 & 30             & 0.114 & 0.144 & 0.120 & 0.088 & 0.091 & 0.080 & 0.062 & 0.069 & 0.050 \\ 
  2 & 2 & 2 & 50             & 0.268 & 0.432 & 0.492 & 0.100 & 0.097 & 0.109 & 0.043 & 0.044 & 0.051 \\ 
  2 & 2 & 6 & 50             & 0.203 & 0.237 & 0.258 & 0.125 & 0.107 & 0.140 & 0.079 & 0.053 & 0.067 \\ 
  3 & 2 & 2 & 20             & 0.134 & 0.158 & 0.173 & 0.066 & 0.074 & 0.062 & 0.045 & 0.060 & 0.042 \\ 
  3 & 2 & 6 & 20             & 0.091 & 0.109 & 0.117 & 0.074 & 0.071 & 0.084 & 0.068 & 0.064 & 0.053 \\ 
  3 & 2 & 2 & 30             & \textbf{0.203} & \textbf{0.330} & \textbf{0.351} & \textbf{0.089} & \textbf{0.088} & \textbf{0.090} & \textbf{0.057} & \textbf{0.044} & \textbf{0.062} \\ 
  3 & 2 & 6 & 30             & 0.164 & 0.208 & 0.186 & 0.102 & 0.113 & 0.072 & 0.072 & 0.063 & 0.056 \\ 
  3 & 2 & 2 & 50             & 0.349 & 0.589 & 0.671 & 0.098 & 0.100 & 0.093 & 0.052 & 0.054 & 0.038 \\ 
  3 & 2 & 6 & 50             & 0.249 & 0.391 & 0.422 & 0.112 & 0.121 & 0.136 & 0.065 & 0.049 & 0.048 \\ 
  5 & 2 & 2 & 20             & 0.166 & 0.243 & 0.288 & 0.084 & 0.072 & 0.071 & 0.045 & 0.051 & 0.056 \\ 
  5 & 2 & 6 & 20             & 0.128 & 0.169 & 0.178 & 0.076 & 0.076 & 0.068 & 0.056 & 0.066 & 0.049 \\ 
  5 & 2 & 2 & 30             & 0.226 & 0.409 & 0.482 & 0.072 & 0.074 & 0.088 & 0.060 & 0.045 & 0.059 \\ 
  5 & 2 & 6 & 30             & 0.199 & 0.298 & 0.323 & 0.091 & 0.099 & 0.101 & 0.052 & 0.060 & 0.063 \\ 
  5 & 2 & 2 & 50             & 0.375 & 0.686 & 0.808 & 0.092 & 0.080 & 0.084 & 0.059 & 0.046 & 0.051 \\ 
  5 & 2 & 6 & 50             & 0.330 & 0.527 & 0.587 & 0.122 & 0.129 & 0.130 & 0.050 & 0.050 & 0.059 \\ 
   \hline
\end{tabular}
\end{table}

\end{document}